\documentclass[twoside]{article}
\usepackage{amsmath}
\usepackage{amssymb}
\usepackage{amsthm}
\usepackage{capt-of}
\usepackage{caption}
\usepackage[strict]{changepage}
\usepackage{chngcntr}
\usepackage[americanvoltage,siunitx]{circuitikz}
\usepackage{color,colortbl}
\usepackage{etoolbox}
\usepackage{fancyhdr}
\usepackage[T1]{fontenc}
\usepackage{gensymb}
\usepackage[margin=1in]{geometry}
\usepackage{graphicx}
\usepackage{hyperref}
\usepackage{import}
\usepackage{indentfirst}
\usepackage{mathptmx}
\usepackage{mathrsfs}
\usepackage{multicol}
\usepackage{multirow}
\usepackage{needspace}
\usepackage{pgfplots}
\usepackage{pgfplotstable}
\usepackage{setspace}
\usepackage{siunitx}
\usepackage{tabu}
\usepackage{tabularx}
\usepackage{tikz}
\usepackage{xspace}

\patchcmd{\thebibliography}{\section*{\refname}}{\vspace{-1em}}{}{}

\usepgfplotslibrary{external}
\usetikzlibrary{positioning,matrix,shapes,chains,arrows}
\newcommand\svgsize[2]{\def\svgwidth{#2}
{\centering\input{#1.pdf_tex}}}
\newcommand\svgc[1]{\def\svgwidth{\columnwidth}
{\centering\input{#1.pdf_tex}}}
\newcommand\svgl[1]{\def\svgwidth{1em}
{\centering\input{#1.pdf_tex}}}
\newcommand\diagrams[0]{\renewcommand\svgsize[2]{\def\svgwidth{##2}
{\centering\input{diagrams/##1.pdf_tex}}}}

\newcommand\pdf[1]{\noindent\includegraphics[width=\columnwidth]{#1.pdf}}
\newcommand\pdfex[1]{\pdf{#1}

\pdf{#1ex}}
\newcommand\pdfmsg[1]{\noindent\begin{minipage}{\columnwidth}\pdf{#1msg}

\pdf{#1}\end{minipage}}
\newcommand\pdfmsgex[1]{\pdfmsg{#1}

\pdf{#1ex}}
\newcommand\code[0]{\renewcommand\pdf[1]{\noindent
\includegraphics[width=\columnwidth]{code/##1.pdf}}}
\newcommand\size[2]{{\fontsize{#1pt}{#1pt}\selectfont#2}}
\newcommand\brokensize[2]{\fontsize{#1pt}{#1pt}\selectfont#2}

\newcounter{paperthmamount}
\newcommand\theorems[0]{
\theoremstyle{remark}
\newtheorem{claim}[subsection]{Claim}
\theoremstyle{plain}
\newtheorem{conjecture}[subsection]{Conjecture}
\theoremstyle{plain}
\newtheorem{corollary}[subsection]{Corollary}
\theoremstyle{definition}
\newtheorem{definition}[subsection]{Definition}
\theoremstyle{plain}
\newtheorem{lemma}[subsection]{Lemma}
\theoremstyle{remark}
\newtheorem{proposition}[subsection]{Proposition}
\theoremstyle{remark}
\newtheorem{remark}[subsection]{Remark}
\theoremstyle{plain}
\newtheorem{theorem}[subsection]{Theorem}
\theoremstyle{definition}
\newtheorem{question}[subsection]{Question}
\newcommand\paperclm[2]
{\begin{claim}\global\expandafter\edef
\csname clm##1\endcsname{Claim \thesubsection\noexpand\xspace}
##2\end{claim}}
\newcommand\papercnj[2]
{\begin{conjecture}\global\expandafter\edef
\csname cnj##1\endcsname{Conjecture \thesubsection\noexpand\xspace}
##2\end{conjecture}}
\newcommand\papercor[2]
{\begin{corollary}\global\expandafter\edef
\csname cor##1\endcsname{Corollary \thesubsection\noexpand\xspace}
##2\end{corollary}}
\newcommand\paperdef[2]
{\begin{definition}\global\expandafter\edef
\csname def##1\endcsname{Definition \thesubsection\noexpand\xspace}
##2\end{definition}}
\newcommand\paperlem[2]
{\begin{lemma}\global\expandafter\edef
\csname lem##1\endcsname{Lemma \thesubsection\noexpand\xspace}
##2\end{lemma}}
\newcommand\paperprp[2]
{\begin{proposition}\global\expandafter\edef
\csname prp##1\endcsname{Proposition \thesubsection\noexpand\xspace}
##2\end{proposition}}
\newcommand\paperqtn[2]
{\begin{question}\global\expandafter\edef
\csname qtn##1\endcsname{Question \thesubsection\noexpand\xspace}
##2\end{question}}
\newcommand\paperrem[2]
{\begin{remark}\global\expandafter\edef
\csname rem##1\endcsname{Remark \thesubsection\noexpand\xspace}
##2\end{remark}}
\newcommand\paperthm[2]
{\begin{theorem}\global\expandafter\edef
\csname thm##1\endcsname{Theorem \thesubsection\noexpand\xspace}
##2\end{theorem}}}
\newcommand\subtheorems[0]{\stepcounter{paperthmamount}
\theoremstyle{remark}
\newtheorem{claim}[subsubsection]{Claim}
\theoremstyle{plain}
\newtheorem{conjecture}[subsubsection]{Conjecture}
\theoremstyle{plain}
\newtheorem{corollary}[subsubsection]{Corollary}
\theoremstyle{definition}
\newtheorem{definition}[subsubsection]{Definition}
\theoremstyle{plain}
\newtheorem{lemma}[subsubsection]{Lemma}
\theoremstyle{remark}
\newtheorem{proposition}[subsubsection]{Proposition}
\theoremstyle{remark}
\newtheorem{remark}[subsubsection]{Remark}
\theoremstyle{plain}
\newtheorem{theorem}[subsubsection]{Theorem}
\theoremstyle{definition}
\newtheorem{question}[subsubsection]{Question}
\newcommand\paperclm[2]
{\begin{claim}\global\expandafter\edef
\csname clm##1\endcsname{Claim \thesubsubsection\noexpand\xspace}
##2\end{claim}}
\newcommand\papercnj[2]
{\begin{conjecture}\global\expandafter\edef
\csname cnj##1\endcsname{Conjecture \thesubsubsection\noexpand\xspace}
##2\end{conjecture}}
\newcommand\papercor[2]
{\begin{corollary}\global\expandafter\edef
\csname cor##1\endcsname{Corollary \thesubsubsection\noexpand\xspace}
##2\end{corollary}}
\newcommand\paperdef[2]
{\begin{definition}\global\expandafter\edef
\csname def##1\endcsname{Definition \thesubsubsection\noexpand\xspace}
##2\end{definition}}
\newcommand\paperlem[2]
{\begin{lemma}\global\expandafter\edef
\csname lem##1\endcsname{Lemma \thesubsubsection\noexpand\xspace}
##2\end{lemma}}
\newcommand\paperprp[2]
{\begin{proposition}\global\expandafter\edef
\csname prp##1\endcsname{Proposition \thesubsubsection\noexpand\xspace}
##2\end{proposition}}
\newcommand\paperqtn[2]
{\begin{question}\global\expandafter\edef
\csname qtn##1\endcsname{Question \thesubsubsection\noexpand\xspace}
##2\end{question}}
\newcommand\paperrem[2]
{\begin{remark}\global\expandafter\edef
\csname rem##1\endcsname{Remark \thesubsubsection\noexpand\xspace}
##2\end{remark}}
\newcommand\paperthm[2]
{\begin{theorem}\global\expandafter\edef
\csname thm##1\endcsname{Theorem \thesubsubsection\noexpand\xspace}
##2\end{theorem}}}

\renewcommand{\headrulewidth}{0pt}
\newcommand\papertitle[1]
{{\centering\fontsize{20pt}{20pt}\textsc{#1}\\\mbox{}\\}
\fancyhead[OC]{\fontsize{12pt}{12pt}\selectfont\textit{#1}}}
\newcounter{people}
\newcommand\paperauthtext[3]{{\centering\fontsize{12pt}{12pt}\selectfont
\textsc{#1}\\[-0.1em]{\fontsize{9pt}{9pt}\selectfont\textit{\ifx&#2&
\vspace{-1em}\else#2\fi}}\\\mbox{}\\
\fancyhead[EC]{\fontsize{12pt}{12pt}\selectfont\textit{#3}}}}
\newcommand\paperauth[2]{{\stepcounter{people}
\ifnum\value{people}=1
{\paperauthtext{#1}{#2}{#1}
\global\def\auth{#1\xspace}}
\else\ifnum\value{people}=2
{\paperauthtext{#1}{#2}{\auth and #1}}
\else{\paperauthtext{#1}{#2}{\auth et al}}\fi\fi}}
\newcommand\physics[0]{
\renewcommand\paperauthtext[4]{{\centering\fontsize{12pt}{12pt}\selectfont
\textsc{##1. ##2}\\[-0.1em]{\fontsize{9pt}{9pt}\selectfont\textit{\ifx&##3&
\vspace{-1em}\else##3\fi}}\\\mbox{}\\
\fancyhead[EC]{\fontsize{12pt}{12pt}\selectfont\textit{##4}}}}
\renewcommand\paperauth[3]{{\stepcounter{people}
\ifnum\value{people}=1
{\paperauthtext{##1}{##2}{##3}{##1. ##2}
\global\def\auth{##2\xspace}}
\else\ifnum\value{people}=2
{\paperauthtext{##1}{##2}{##3}{\auth and ##2}}
\else{\paperauthtext{##1}{##2}{##3}{\auth et al}}\fi\fi}}}
\newcommand\paperdate[1]{{\centering\fontsize{9pt}{9pt}\selectfont\text{
(Received #1)}\\[2em]}}

\newcommand{\paperhead}[1]{\fancyhead[EC]{\fontsize{12pt}{12pt}\selectfont
\textit{#1}}}
\makeatletter
\newenvironment{paperadjustwidth}[2]{
  \begin{list}{}{
    \setlength\partopsep\z@
    \setlength\topsep\z@
    \setlength\listparindent\parindent
    \setlength\parsep\parskip
    \linespread{0.75}\selectfont
    \@ifmtarg{#1}{\setlength{\leftmargin}{\z@}}
                 {\setlength{\leftmargin}{#1}}
    \@ifmtarg{#2}{\setlength{\rightmargin}{\z@}}
                 {\setlength{\rightmargin}{#2}}
    }
    \item[]}{\end{list}}
\makeatother

\newenvironment{paperabs}
{\begin{paperadjustwidth}{0.5in}{0.5in}\bgroup\fontsize{9pt}{9pt}\selectfont
\hspace{0.5in}}
{\egroup\end{paperadjustwidth}}

\newenvironment{paper}
{\begin{multicols*}{2}\bgroup\fontsize{12pt}{12pt}\selectfont}
{\egroup\end{multicols*}}

\newsavebox{\sourcebox}
\newcommand{\papersource}[1]{
\vspace{-2em}
\text{}\\*
\fontsize{9pt}{9pt}\selectfont
\noindent\renewcommand{\labelenumi}{}
\savebox{\sourcebox}{\parbox{3in}{\begin{enumerate}
\setlength{\leftmargini}{-1ex}
\setlength{\leftmargin}{-1ex}
\setlength{\labelwidth}{0pt}
\setlength{\labelsep}{0pt}
\setlength{\listparindent}{0pt}
\item\textit{\hspace{-0.35in}#1}
\end{enumerate}}}
\usebox{\sourcebox}
}

\newcounter{paperseccounter}
\newcounter{papersubseccounter}[paperseccounter]
\newcommand\papersec[1]{\needspace{1in}
\stepcounter{paperseccounter}
\stepcounter{section}
\begin{center}\Roman{paperseccounter} \textsc{#1}\end{center}}
\newcommand\papersubsec[1]{\needspace{1in}
\stepcounter{papersubseccounter}
\addtocounter{subsection}{\thepaperthmamount}
\setcounter{subsubsection}{0}
{\begin{center}
\Roman{section}.\Roman{papersubseccounter}
\textsc{#1}\\[0.5em]\end{center}}}

\newcounter{papereqcounter}
\newcommand\papereq[3]{{
\stepcounter{papereqcounter}
\mbox{}\vspace{-0.75em}
\begin{equation*}
#2
\tag*{\fontsize{12pt}{12pt}\selectfont
$\begin{array}{r}
\cr{\text{(\arabic{papereqcounter})}}
\cr{\fontsize{9pt}{9pt}\selectfont\textit{\ifx\\#3\\~\else(\fi#3\ifx\\#3\\~
\else)\fi}}
\end{array}$}
\end{equation*}
}
\expandafter\edef\csname eq#1\endcsname{(\arabic{papereqcounter})\noexpand
\xspace}}

\newcommand{\papervar}[3]
{&$#1$ & #2 \ifx\\#3\\~\else($\smash{\text{\si{\fi
#3\ifx\\#3\\~\else}}}$)\fi\\}
\newenvironment{paperwhere}
{\begin{minipage}{\columnwidth}
\bgroup\fontsize{9pt}{9pt}\selectfont Where:\vspace{2pt}\\\begin{tabular}
{rr@{ = }p{\linewidth}}}
{\end{tabular}\egroup\end{minipage}\vspace{5pt}}

\definecolor{LineGray}{gray}{0.5}
\newtabulinestyle{outer=2.25pt LineGray}
\newtabulinestyle{inner=0.75pt LineGray}
\newcommand{\paperiline}[0]{\tabucline[inner]{-}}
\newcommand{\paperoline}[0]{\tabucline[outer]{-}}

\newcolumntype{I}{X[-5,c]}
\newcolumntype{U}{@{}X[-5,r]@{$\pm$}X[-5,l]@{}}
\newcolumntype{C}{@{}X[-5,c]@{}}

\newcounter{papertableindexcounter}
\newcommand{\papertableindexheader}[0]{\multirow{2}{*}{\textsc{Index}}}
\newcommand{\papertableindex}[0]{\stepcounter{papertableindexcounter}
\arabic{papertableindexcounter}}
\newcommand{\papertableuheadersymbol}[1]{&\multicolumn{2}{c|[inner]}{$#1$}}
\newcommand{\papertableuheadersymbole}[1]{&\multicolumn{2}{c|[outer]}{$#1$}}
\newcommand{\papertableuheaderunit}[1]{&\multicolumn{2}{c|[inner]}{(#1)}}
\newcommand{\papertableuheaderunite}[1]{&\multicolumn{2}{c|[outer]}{(#1)}}
\newcommand{\papertablecheadersymbol}[1]{&$#1$}
\newcommand{\papertablecheaderunit}[2]{&($\pm$#1 #2)}

\newcommand{\papertableuval}[2]{& #1 & #2}
\newcommand{\papertablecval}[1]{& #1}

\newcounter{paperfigurecounter}
\newcommand{\papercap}[2]{\bgroup\stepcounter{paperfigurecounter}
\captionof{figure}{\fontsize{9pt}{9pt}\selectfont
\hspace{0.3in}Fig.~\arabic{paperfigurecounter}.\quad#2}
\egroup\expandafter\edef
\csname fig#1\endcsname{Fig.~\arabic{paperfigurecounter}\noexpand\xspace}}

\newcommand\paperfig[3]{\noindent\begin{minipage}{\columnwidth}
#2\papercap{#1}{#3}\end{minipage}\expandafter\edef
\csname fig#1\endcsname{Fig.~\arabic{paperfigurecounter}\noexpand\xspace}}
\newcommand\papersvg[3]{\paperfig{#1}{\svgc{#2}}{#3}}

\newcommand{\paperaxis}[9]
{title=#1,
axis x line = bottom,
xmin=#4,xmax=#6,
axis y line = left,
ymin=#5,ymax=#7,
height = 180pt,
grid=both,
x axis line style=-,
y axis line style=-,
x tick label style={
/pgf/number format/.cd,
fixed,
fixed zerofill,
precision=#8,
/tikz/.cd},
y tick label style={
/pgf/number format/.cd,
fixed,
fixed zerofill,
precision=#9,
/tikz/.cd}}
\newcommand{\paperaxisxlabel}[2]{
xlabel=\fontsize{10pt}{10pt}\selectfont#1$(#2)\rightarrow$}
\newcommand{\paperaxisylabel}[2]{
ylabel=\fontsize{10pt}{10pt}\selectfont#1$(#2)\rightarrow$}
\newcommand{\papergraphoutline}[4]{
\addplot [mark=none,line width=0.75pt] coordinates {
(#1,#2)
(#1,#4)
(#3,#4)
(#3,#2)
(#1,#2)};}

\newcommand{\comment}[1]{}

\newcommand{\abs}[1]{\left\lvert#1\right\rvert}
\newcommand{\oo}[0]{\infty}
\newcommand{\sigmaSum}[3]{\sum\limits_{#1}^{#2} #3}
\newcommand{\limto}[3]{\lim\limits_{#1\rightarrow#2}#3}
\renewcommand{\d}[0]{\mathrm{d}}
\newcommand{\cross}[0]{\times}
\newcommand{\lp}{\left(}
\newcommand{\rp}{\right)}
\newcommand\pars[1]{\lp#1\rp}
\newcommand\sqbrack[1]{\left[#1\right]}
\newcommand\R{\mathbb{R}}
\newcommand\di{\partial}
\newcommand\x{\times}
\newcommand\del{\nabla}
\code
\renewcommand\svgsize[2]{\def\svgwidth{##2}
{\centering\input{diagrams/##1.pdf_tex}}}
\begin{document}
\papertitle{The 250 Knots with up to 10 Crossings}
\paperauth{Andrey Boris Khesin}{University of Toronto}
\begin{paperabs}
The list of knots with up to 10 crossings is commonly referred to as the Rolfsen
Table.
This paper presents a way to generate the Rolfsen table in a simple, clear, and
reproducible manner.
The methods we use are similar to those used by J.~Hoste, M.~Thistlethwaite, and
J.~Weeks in \cite{htw}.
The difference between our methods comes from the fact that \cite{htw} uses a
more complicated algorithm to be able to find all the knots with up to 17
crossings, while our approach demonstrates a simpler way to find the knots up to
10 crossings.
We do this by generating all planar knot diagrams with up to 10 crossings and
applying several simplifications to group the knot diagrams into equivalence
classes.
From these classes, we generate the full list of candidate knots and reduce it
with several sets of moves.
Lastly, we use invariants to show that each of the 250 diagrams generated is
distinct,  proving that there are exactly 250 knots with 10 crossings or fewer.
Though the algorithms used could be made more efficient, readability was chosen
over speed for simplicity and reproducibility.
\end{paperabs}
\begin{paper}
\papersec{Introduction}

The Rolfsen table is the list of the 250 knots with 10 crossings or fewer.
Here we attempt to generate it and prove its completeness by using a computer
algorithm.
This has been accomplished several times in the past.
A notable example is \cite{htw}, where J.~Hoste, M.~Thistlethwaite, and J.~Weeks
found all of the knots with up to 17 crossings.
The methods we use are far less advanced, which allows us to effect a less
intensive computation, finding the knots with up to 10 crossings, but use a
simpler algorithm to do accomplish this.

Although it is possible to compute the Rolfsen table by hand, it is a rather
tedious task.
Our calculation is made possible by using a computer.
To demonstrate a method of generating the Rolfsen table, we create a simple
algorithm for finding all 250 knots with up to 10 crossings, partially
sacrificing efficiency in the process.

We begin our reconstruction of the Rolfsen table by considering which knot
diagrams could potentially be included in the table.
There are only a finite number of ways to draw a knot diagram with a given
number of crossings.
Additionally, many of these knot diagrams are \textit{reducible}, which means
that they are equivalent to other knot diagrams with fewer crossings.

There are far more than 250 knot diagrams with up to 10 crossings, even after
only irreducible knot diagrams are considered.
The reason for this is that there are several moves that can transform one knot
diagram into an equivalent one.
Two knot diagrams are \textit{equivalent} if and only if there exists a series
of such moves that transforms one of the diagrams into the other.

\begin{center}\begin{minipage}{\columnwidth}\begin{center}
Reidemeister Moves\vspace{0.5em}
\def\svgwidth{0.65\columnwidth}
{\centering
\begingroup%
  \makeatletter%
  \providecommand\color[2][]{%
    \errmessage{(Inkscape) Color is used for the text in Inkscape, but the package 'color.sty' is not loaded}%
    \renewcommand\color[2][]{}%
  }%
  \providecommand\transparent[1]{%
    \errmessage{(Inkscape) Transparency is used (non-zero) for the text in Inkscape, but the package 'transparent.sty' is not loaded}%
    \renewcommand\transparent[1]{}%
  }%
  \providecommand\rotatebox[2]{#2}%
  \ifx\svgwidth\undefined%
    \setlength{\unitlength}{493.64649833bp}%
    \ifx\svgscale\undefined%
      \relax%
    \else%
      \setlength{\unitlength}{\unitlength * \real{\svgscale}}%
    \fi%
  \else%
    \setlength{\unitlength}{\svgwidth}%
  \fi%
  \global\let\svgwidth\undefined%
  \global\let\svgscale\undefined%
  \makeatother%
  \begin{picture}(1,0.52383986)%
    \put(0,0){\includegraphics[width=\unitlength,page=1]{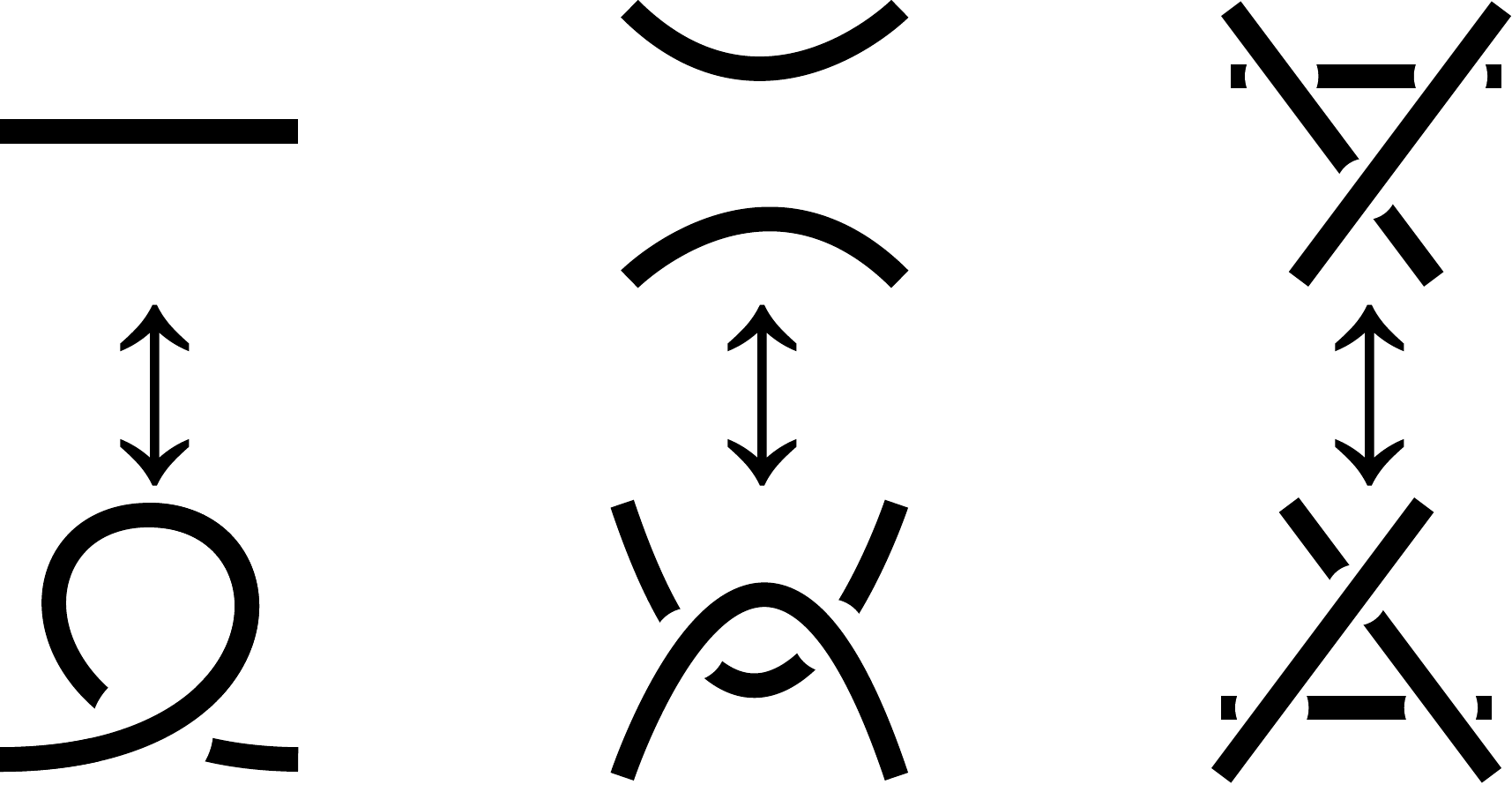}}%
  \end{picture}%
\endgroup%
}\\
First\hspace{0.14\columnwidth}Second\hspace{0.11\columnwidth}Third\\
\end{center}\end{minipage}\end{center}
\begin{center}\begin{minipage}{\columnwidth}\begin{center}
Crossing Number-Preserving Moves\vspace{0.5em}
\def\svgwidth{0.65\columnwidth}
{\centering
\begingroup%
  \makeatletter%
  \providecommand\color[2][]{%
    \errmessage{(Inkscape) Color is used for the text in Inkscape, but the package 'color.sty' is not loaded}%
    \renewcommand\color[2][]{}%
  }%
  \providecommand\transparent[1]{%
    \errmessage{(Inkscape) Transparency is used (non-zero) for the text in Inkscape, but the package 'transparent.sty' is not loaded}%
    \renewcommand\transparent[1]{}%
  }%
  \providecommand\rotatebox[2]{#2}%
  \ifx\svgwidth\undefined%
    \setlength{\unitlength}{529.74258344bp}%
    \ifx\svgscale\undefined%
      \relax%
    \else%
      \setlength{\unitlength}{\unitlength * \real{\svgscale}}%
    \fi%
  \else%
    \setlength{\unitlength}{\svgwidth}%
  \fi%
  \global\let\svgwidth\undefined%
  \global\let\svgscale\undefined%
  \makeatother%
  \begin{picture}(1,0.80701776)%
    \put(0,0){\includegraphics[width=\unitlength,page=1]{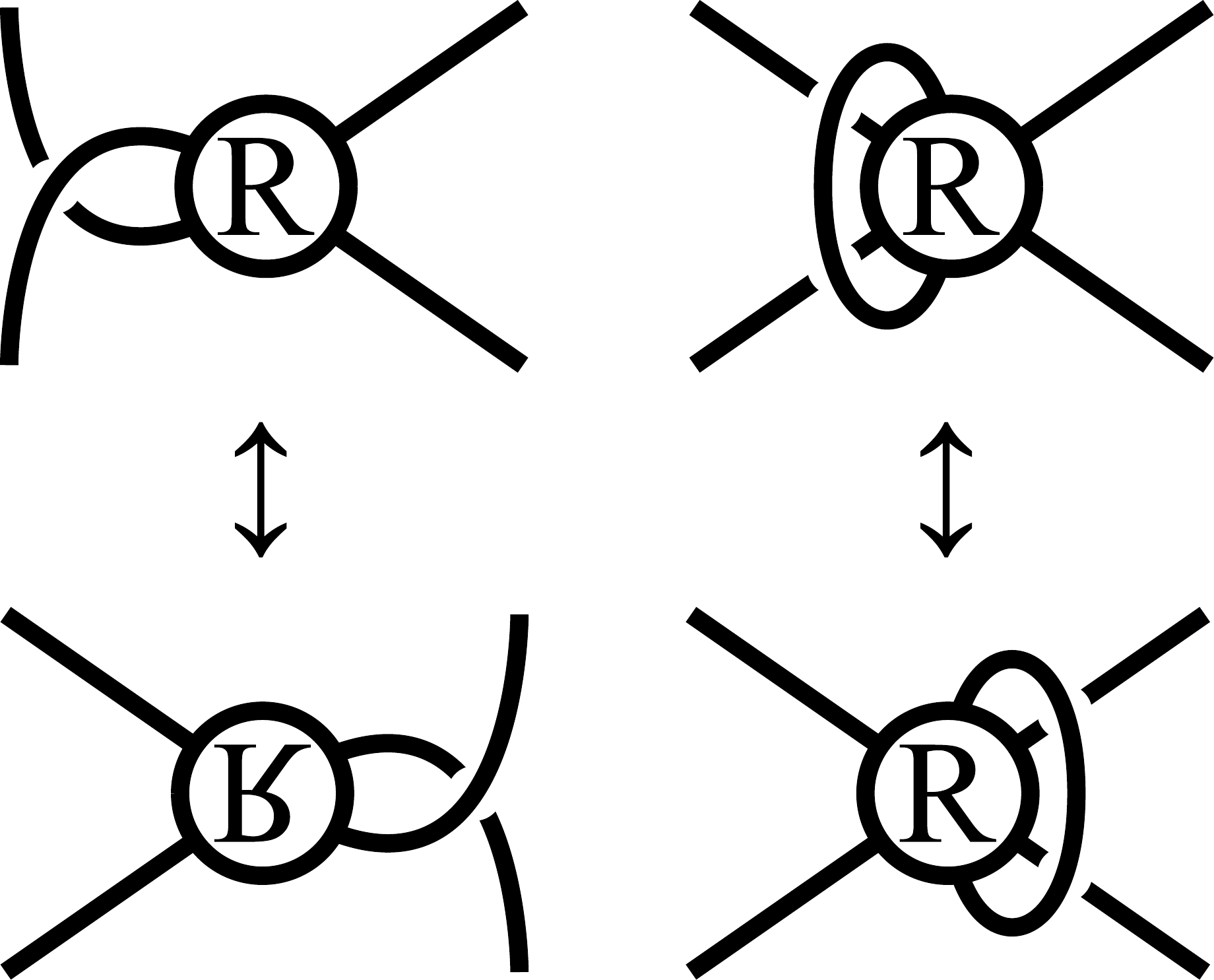}}%
  \end{picture}%
\endgroup%
}\\
Flype\hspace{0.225\columnwidth}2--Pass\\
\end{center}\end{minipage}\end{center}
\begin{center}\begin{minipage}{\columnwidth}\begin{center}
Crossing Number-Reducing Moves\vspace{0.5em}
\def\svgwidth{0.65\columnwidth}
{\centering
\begingroup%
  \makeatletter%
  \providecommand\color[2][]{%
    \errmessage{(Inkscape) Color is used for the text in Inkscape, but the package 'color.sty' is not loaded}%
    \renewcommand\color[2][]{}%
  }%
  \providecommand\transparent[1]{%
    \errmessage{(Inkscape) Transparency is used (non-zero) for the text in Inkscape, but the package 'transparent.sty' is not loaded}%
    \renewcommand\transparent[1]{}%
  }%
  \providecommand\rotatebox[2]{#2}%
  \ifx\svgwidth\undefined%
    \setlength{\unitlength}{527.26144926bp}%
    \ifx\svgscale\undefined%
      \relax%
    \else%
      \setlength{\unitlength}{\unitlength * \real{\svgscale}}%
    \fi%
  \else%
    \setlength{\unitlength}{\svgwidth}%
  \fi%
  \global\let\svgwidth\undefined%
  \global\let\svgscale\undefined%
  \makeatother%
  \begin{picture}(1,0.81081534)%
    \put(0,0){\includegraphics[width=\unitlength,page=1]{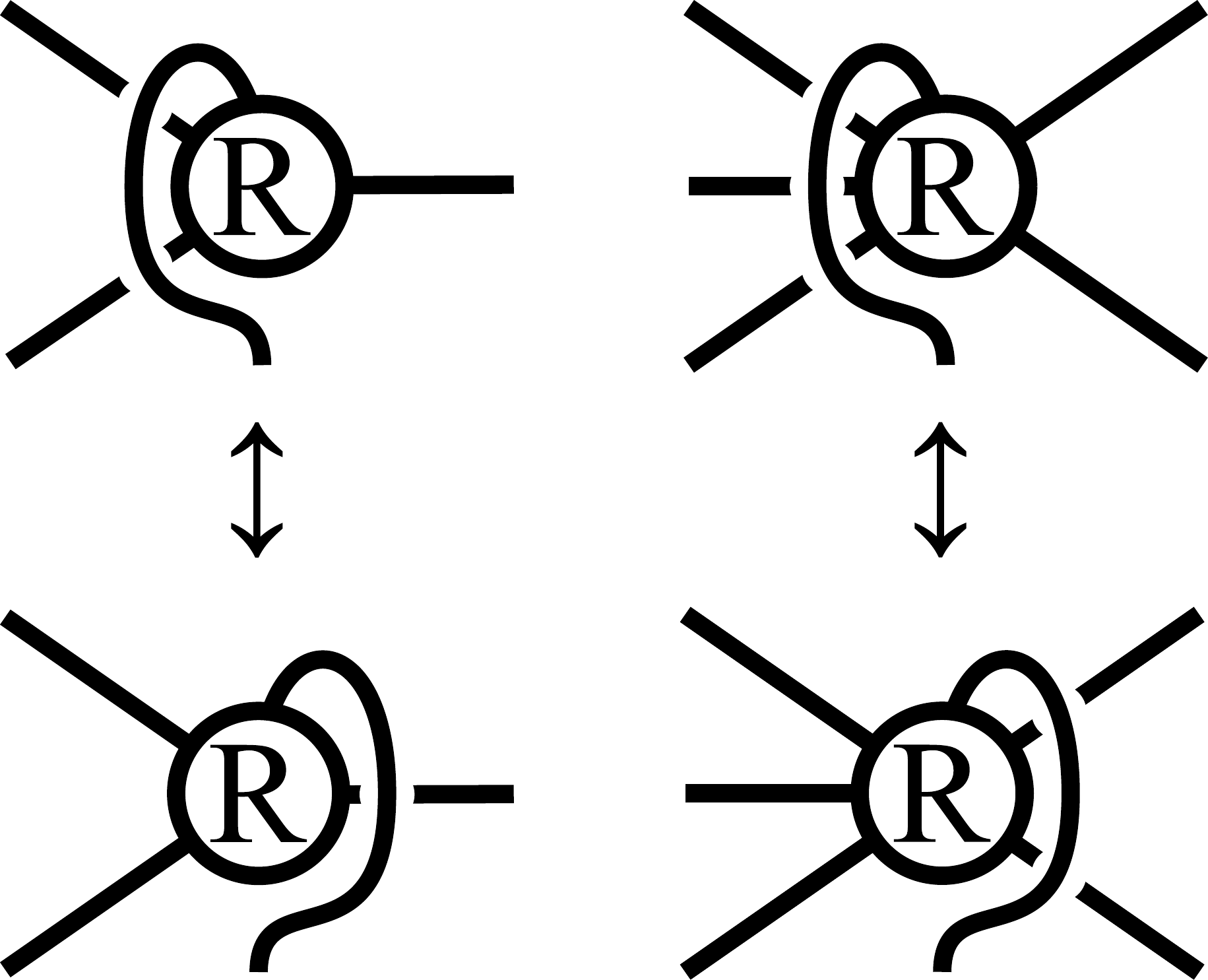}}%
  \end{picture}%
\endgroup%
}\\
(2, 1)--Pass\hspace{0.125\columnwidth}(3, 2)--Pass
\end{center}\end{minipage}\end{center}

\paperfig{Moves}{}
{The 6 moves that we use to construct the Rolfsen table, as well as the second
Reidemeister move.
The letter R is used to denote a tangle with an appropriate number of strands.
If the letter R appears in a different orientation it is because the move caused
the corresponding part of the knot diagram to flip.}

\paperfig{Trefoil}{
\def\svgwidth{0.4\columnwidth}
{\centering
\begingroup%
  \makeatletter%
  \providecommand\color[2][]{%
    \errmessage{(Inkscape) Color is used for the text in Inkscape, but the package 'color.sty' is not loaded}%
    \renewcommand\color[2][]{}%
  }%
  \providecommand\transparent[1]{%
    \errmessage{(Inkscape) Transparency is used (non-zero) for the text in Inkscape, but the package 'transparent.sty' is not loaded}%
    \renewcommand\transparent[1]{}%
  }%
  \providecommand\rotatebox[2]{#2}%
  \ifx\svgwidth\undefined%
    \setlength{\unitlength}{496.72726665bp}%
    \ifx\svgscale\undefined%
      \relax%
    \else%
      \setlength{\unitlength}{\unitlength * \real{\svgscale}}%
    \fi%
  \else%
    \setlength{\unitlength}{\svgwidth}%
  \fi%
  \global\let\svgwidth\undefined%
  \global\let\svgscale\undefined%
  \makeatother%
  \begin{picture}(1,0.94159925)%
    \put(0,0){\includegraphics[width=\unitlength,page=1]{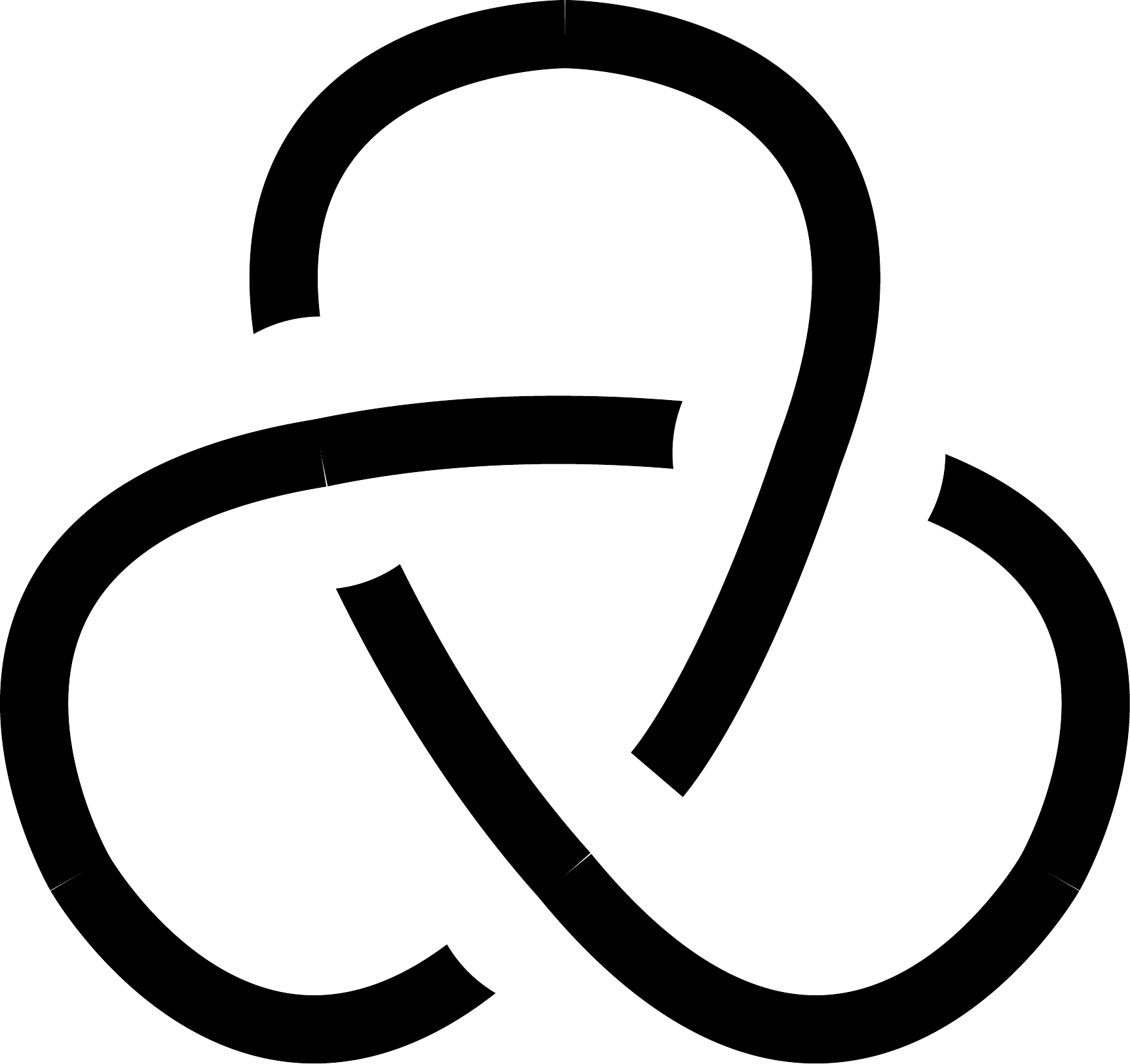}}%
  \end{picture}%
\endgroup%
}\hfill\def\svgwidth{0.4\columnwidth}
{\centering
\begingroup%
  \makeatletter%
  \providecommand\color[2][]{%
    \errmessage{(Inkscape) Color is used for the text in Inkscape, but the package 'color.sty' is not loaded}%
    \renewcommand\color[2][]{}%
  }%
  \providecommand\transparent[1]{%
    \errmessage{(Inkscape) Transparency is used (non-zero) for the text in Inkscape, but the package 'transparent.sty' is not loaded}%
    \renewcommand\transparent[1]{}%
  }%
  \providecommand\rotatebox[2]{#2}%
  \ifx\svgwidth\undefined%
    \setlength{\unitlength}{496.72722813bp}%
    \ifx\svgscale\undefined%
      \relax%
    \else%
      \setlength{\unitlength}{\unitlength * \real{\svgscale}}%
    \fi%
  \else%
    \setlength{\unitlength}{\svgwidth}%
  \fi%
  \global\let\svgwidth\undefined%
  \global\let\svgscale\undefined%
  \makeatother%
  \begin{picture}(1,0.9415993)%
    \put(0,0){\includegraphics[width=\unitlength,page=1]{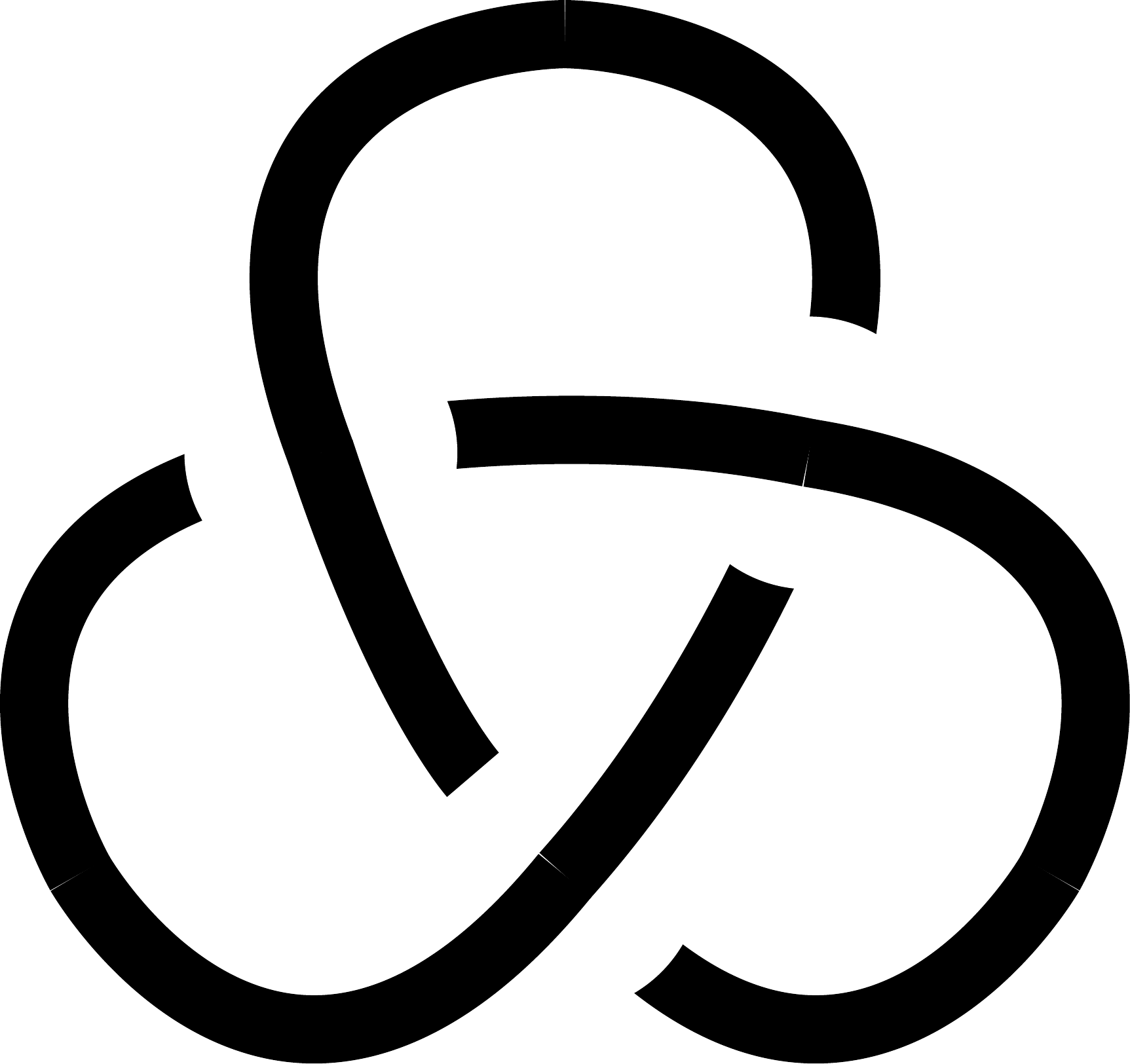}}%
  \end{picture}%
\endgroup%
}
\begin{center}Right\hspace{0.5\columnwidth}Left\end{center}}
{The right-handed trefoil and the left-handed trefoil.
These knots are considered equivalent for our purposes as they are mirror images
of each other.
However, it is important to note that no series of moves can transform one of
these into the other, so while they are not equivalent knots, we only include
one of them in the Rolfsen table.}

The first manner in which we simplify the list of knot diagrams is by
eliminating knot diagrams that are mirror images of each other.
For example, the right-handed and left-handed trefoils are not equivalent as it
is impossible to turn one into the other (see \figTrefoil).
We only include one of the two in the Rolfsen table.
The notation we use to represent a knot diagram does not encode the handedness
of the knot so this is not an issue.\\

\paperfig{Composite}
{\begin{center}\def\svgwidth{0.7\columnwidth}
{\centering
\begingroup%
  \makeatletter%
  \providecommand\color[2][]{%
    \errmessage{(Inkscape) Color is used for the text in Inkscape, but the package 'color.sty' is not loaded}%
    \renewcommand\color[2][]{}%
  }%
  \providecommand\transparent[1]{%
    \errmessage{(Inkscape) Transparency is used (non-zero) for the text in Inkscape, but the package 'transparent.sty' is not loaded}%
    \renewcommand\transparent[1]{}%
  }%
  \providecommand\rotatebox[2]{#2}%
  \ifx\svgwidth\undefined%
    \setlength{\unitlength}{497.63775732bp}%
    \ifx\svgscale\undefined%
      \relax%
    \else%
      \setlength{\unitlength}{\unitlength * \real{\svgscale}}%
    \fi%
  \else%
    \setlength{\unitlength}{\svgwidth}%
  \fi%
  \global\let\svgwidth\undefined%
  \global\let\svgscale\undefined%
  \makeatother%
  \begin{picture}(1,0.84588607)%
    \put(0,0){\includegraphics[width=\unitlength,page=1]{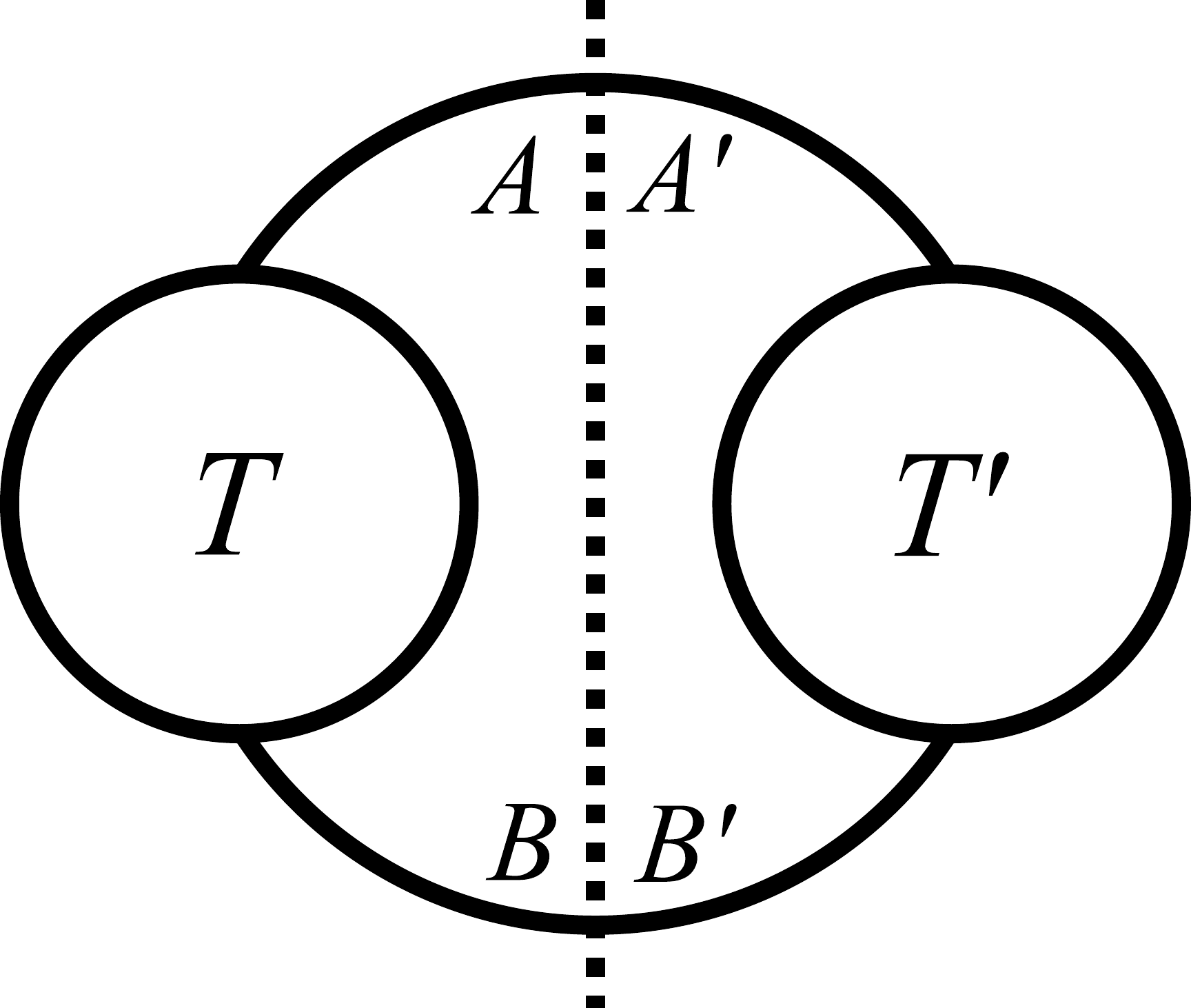}}%
  \end{picture}%
\endgroup%
}\end{center}\vspace{-1em}}
{An example of a composite knot diagram.
This knot diagram can be cut along the dotted line into two knot factors, $T$
and $T'$.
Both $T$ and $T'$ can be cut along one of their edges to create the pairs of
ends $A$ and $B$, as well as $A'$ and $B'$, respectively.
If $A$ is joined to $A'$ and $B$ to $B'$, the resulting knot is the knot
composition of $T$ and $T'$, which are its knot factors.
Knot diagrams that cannot be decomposed into two such knot factors are prime and
are the kind of knot diagrams that we want to include in our tabulation.}

For any two knot diagrams $T$ and $T'$, we cut $T$ at some point to create ends
$A$ and $B$ and we cut $T'$ at some point to create ends $A'$ and $B'$.
Joining $A$ to $A'$ and $B$ to $B'$ results in one larger knot, $R$.
We say that $T$ and $T'$ are \textit{knot factors} of $R$.
The commutative operation of joining $T$ and $T'$ to create $R$ is called
\textit{knot composition}.
Knots that cannot be decomposed into two knot factors other than themselves and
the unknot are called \textit{prime}.
If they can be decomposed this way, they are called \textit{composite}
(see \figComposite).
We only include prime knots in our tabulation.

Lastly, we group knots into equivalence classes based on whether or not there
exists a series of moves that transform one knot diagram into another (see
\figMoves).
Out of each of these equivalence classes, we select one knot to include in our
tabulation.
Having done this, all that remains is proving that our list contains no
remaining equivalent knot diagrams.

\papersec{MD Codes}

The first thing we need to do is to establish is a way to efficiently represent
a knot diagram with some sort of notation.
This notation must be relatively easy for both humans and computers to work
with.

In \cite{htw}, a notation is used to represent an $n$-crossing knot diagram with
$n$ integers.
This notation is called Dowker notation.
Its density and simplicity make it convenient for our purposes.

For an $n$-crossing knot diagram, its representation in Dowker notation is
constructed as follows.
We start by picking an arbitrary point on one of the knot diagram's edges as
well as an arbitrary direction along that edge.
We then move along the knot diagram, moving along each edge, until we have
traveled along all $2n$ edges and have returned to our starting point.
Note that we will pass each crossing twice, once under and once over.
Each time we pass a crossing, we consider the number of crossings that we
have encountered so far and write that number down at the crossing that we are
passing.
In other words, when we encounter our first crossing, we write down the number 1
at that point.
It is important to distinguish between writing the number on the upper or the
lower strand of a crossing.
If we passed the first crossing while traveling along the upper strand, we write
down the number 1 on the upper strand and vice versa.
Continuing, we would write down the number 2 at the next crossing we encounter.
We would end up writing each number from 1 to $2n$ exactly once.
Furthermore, these numbers would be grouped into $n$ pairs, as there would be
two numbers written at each of the $n$ crossings.

Note that since any two closed curves intersect in an even number of places, it
follows that the pair of numbers written at each crossing will contain one odd
number and one even number.
If this were not the case then we would be able to leave a crossing, travel in a
closed loop, and come back to that crossing having encountered an odd number of
crossings along the way, which is not possible.

The $n$ pairs of numbers have no order, so sorting them in ascending order by
comparing the odd value in each pair does not sacrifice any information.
It then follows that the list of even values, sorted by their corresponding odd
value, is sufficient to fully reconstruct the original list of pairs.\\

\paperfig{Labeled}
{\begin{center}\def\svgwidth{0.5\columnwidth}
{\centering
\begingroup%
  \makeatletter%
  \providecommand\color[2][]{%
    \errmessage{(Inkscape) Color is used for the text in Inkscape, but the package 'color.sty' is not loaded}%
    \renewcommand\color[2][]{}%
  }%
  \providecommand\transparent[1]{%
    \errmessage{(Inkscape) Transparency is used (non-zero) for the text in Inkscape, but the package 'transparent.sty' is not loaded}%
    \renewcommand\transparent[1]{}%
  }%
  \providecommand\rotatebox[2]{#2}%
  \ifx\svgwidth\undefined%
    \setlength{\unitlength}{496.72724986bp}%
    \ifx\svgscale\undefined%
      \relax%
    \else%
      \setlength{\unitlength}{\unitlength * \real{\svgscale}}%
    \fi%
  \else%
    \setlength{\unitlength}{\svgwidth}%
  \fi%
  \global\let\svgwidth\undefined%
  \global\let\svgscale\undefined%
  \makeatother%
  \begin{picture}(1,0.95164655)%
    \put(0,0){\includegraphics[width=\unitlength,page=1]{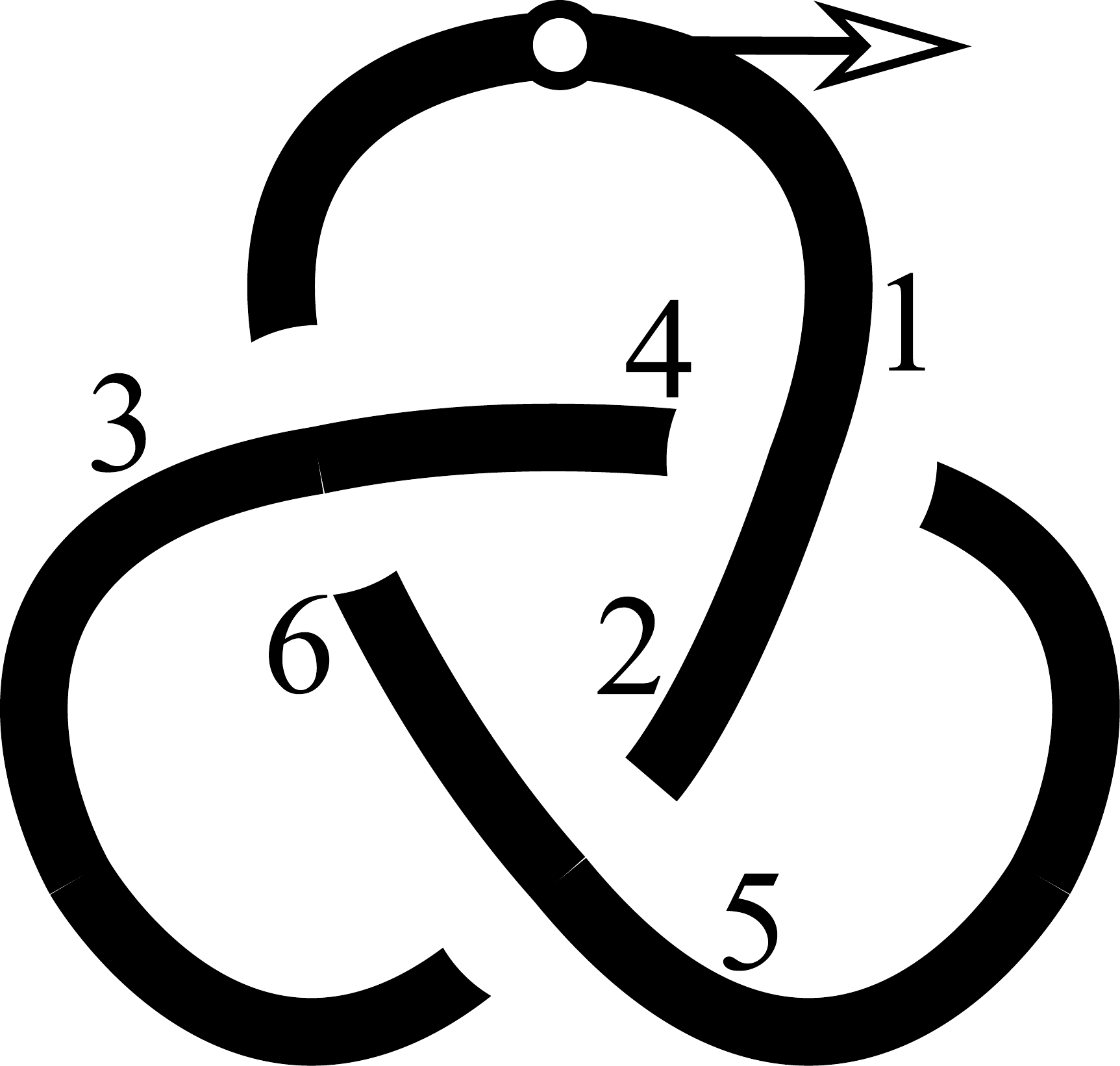}}%
  \end{picture}%
\endgroup%
}\end{center}}
{The right-handed trefoil with the strands in its crossings labeled from 1 to 6.
The labeling starts at the white circle in the centre of the top edge and
proceeds to the right in the direction of the arrow.
The labeling continues until all 6 strands at the knot diagram's crossings are
labeled.
We see that the pairs, when sorted in ascending order by their odd value, are
(1, 4), (3, 6), and (5, 2).}

As an example, we show how we would find the Dowker notation for the trefoil.
Note that the handedness of the trefoil is irrelevant.
After labeling the trefoil, the pairs are (1, 4), (2, 5) and (3, 6) (see
\figLabeled).
We reorder the values in some of the pairs, in this case in (2, 5), to place the
odd value first.
Then the pairs can be ordered by their odd value to get (1, 4), (3, 6), and
(5, 2).
The original pairs can be reconstructed with the sequence (4, 6, 2) as there is
a unique way of reestablishing the odd counterparts to each of the even numbers.
Thus, (4, 6, 2) is the Dowker notation for the trefoil.
Since this sequence contains only even numbers, storing half of each value works
just as well and makes some computations easier.
Therefore, we represent the trefoil by (2, 3, 1).
We call the notation that stores half of each integer an MD code (M is for
modified).
The $2\times n$ matrix of pairs is called an ED code (E is for extended).
The ED code for the trefoil is $\begin{pmatrix}1&3&5\\4&6&2\end{pmatrix}$.

We will later refer to examining permutations.
Since a MD code is a permutation of the numbers from 1 to $n$, we can examine
each such permutation to see if it encodes a viable knot diagram.
Note the distinction between a permutation, one of the many ways of ordering the
numbers from 1 to $n$, and an MD code, a permutation of the number from 1 to $n$
that encodes a particular knot diagram.\\

\paperfig{Crossings}
{\begin{center}Crossings\end{center}\vspace{-1em}
\def\svgwidth{0.33\columnwidth}
{\centering
\begingroup%
  \makeatletter%
  \providecommand\color[2][]{%
    \errmessage{(Inkscape) Color is used for the text in Inkscape, but the package 'color.sty' is not loaded}%
    \renewcommand\color[2][]{}%
  }%
  \providecommand\transparent[1]{%
    \errmessage{(Inkscape) Transparency is used (non-zero) for the text in Inkscape, but the package 'transparent.sty' is not loaded}%
    \renewcommand\transparent[1]{}%
  }%
  \providecommand\rotatebox[2]{#2}%
  \ifx\svgwidth\undefined%
    \setlength{\unitlength}{316.54126374bp}%
    \ifx\svgscale\undefined%
      \relax%
    \else%
      \setlength{\unitlength}{\unitlength * \real{\svgscale}}%
    \fi%
  \else%
    \setlength{\unitlength}{\svgwidth}%
  \fi%
  \global\let\svgwidth\undefined%
  \global\let\svgscale\undefined%
  \makeatother%
  \begin{picture}(1,0.98976115)%
    \put(0,0){\includegraphics[width=\unitlength,page=1]{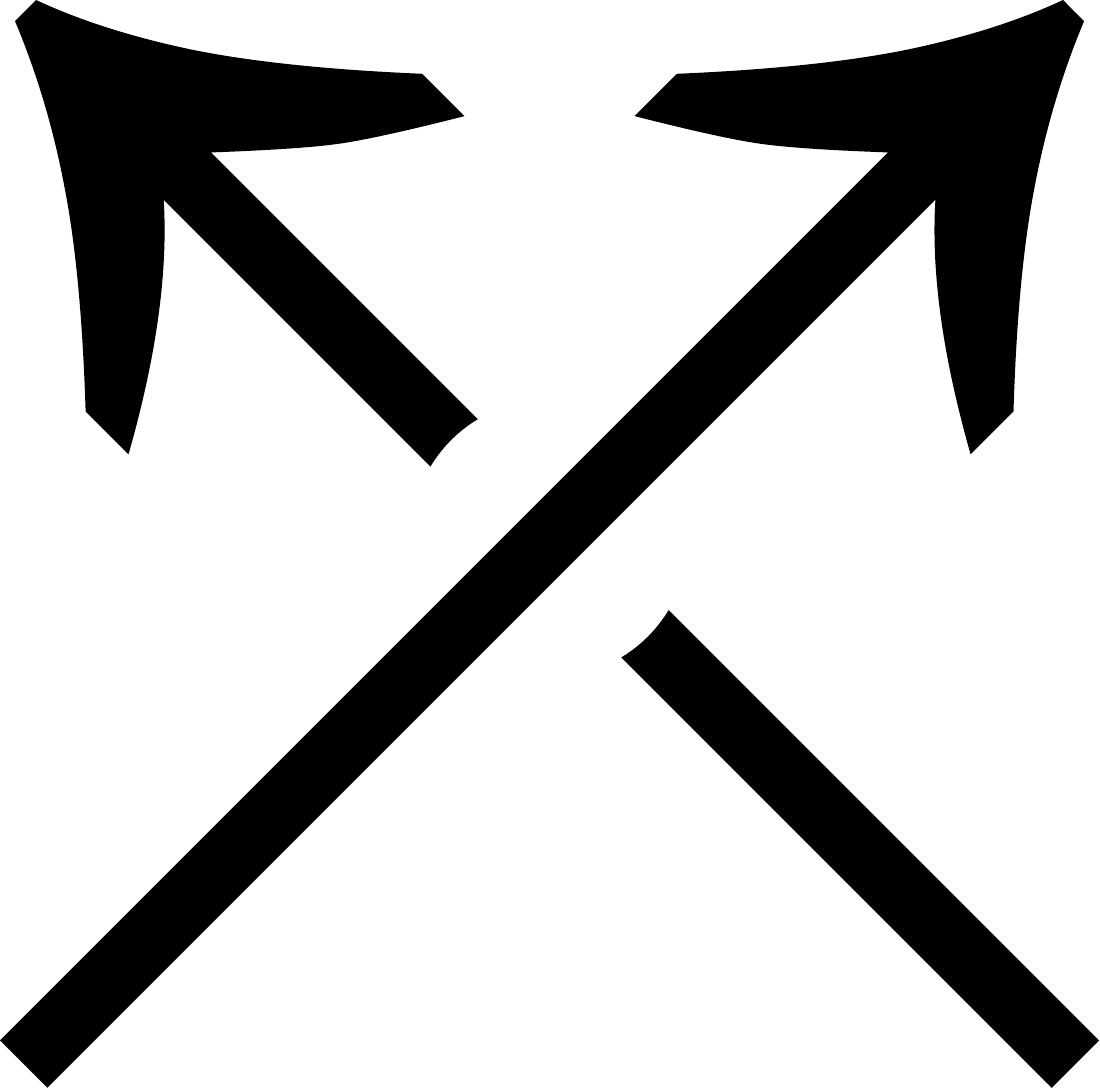}}%
  \end{picture}%
\endgroup%
}
\hfill
\def\svgwidth{0.33\columnwidth}
{\centering
\begingroup%
  \makeatletter%
  \providecommand\color[2][]{%
    \errmessage{(Inkscape) Color is used for the text in Inkscape, but the package 'color.sty' is not loaded}%
    \renewcommand\color[2][]{}%
  }%
  \providecommand\transparent[1]{%
    \errmessage{(Inkscape) Transparency is used (non-zero) for the text in Inkscape, but the package 'transparent.sty' is not loaded}%
    \renewcommand\transparent[1]{}%
  }%
  \providecommand\rotatebox[2]{#2}%
  \ifx\svgwidth\undefined%
    \setlength{\unitlength}{316.54130073bp}%
    \ifx\svgscale\undefined%
      \relax%
    \else%
      \setlength{\unitlength}{\unitlength * \real{\svgscale}}%
    \fi%
  \else%
    \setlength{\unitlength}{\svgwidth}%
  \fi%
  \global\let\svgwidth\undefined%
  \global\let\svgscale\undefined%
  \makeatother%
  \begin{picture}(1,0.9897609)%
    \put(0,0){\includegraphics[width=\unitlength,page=1]{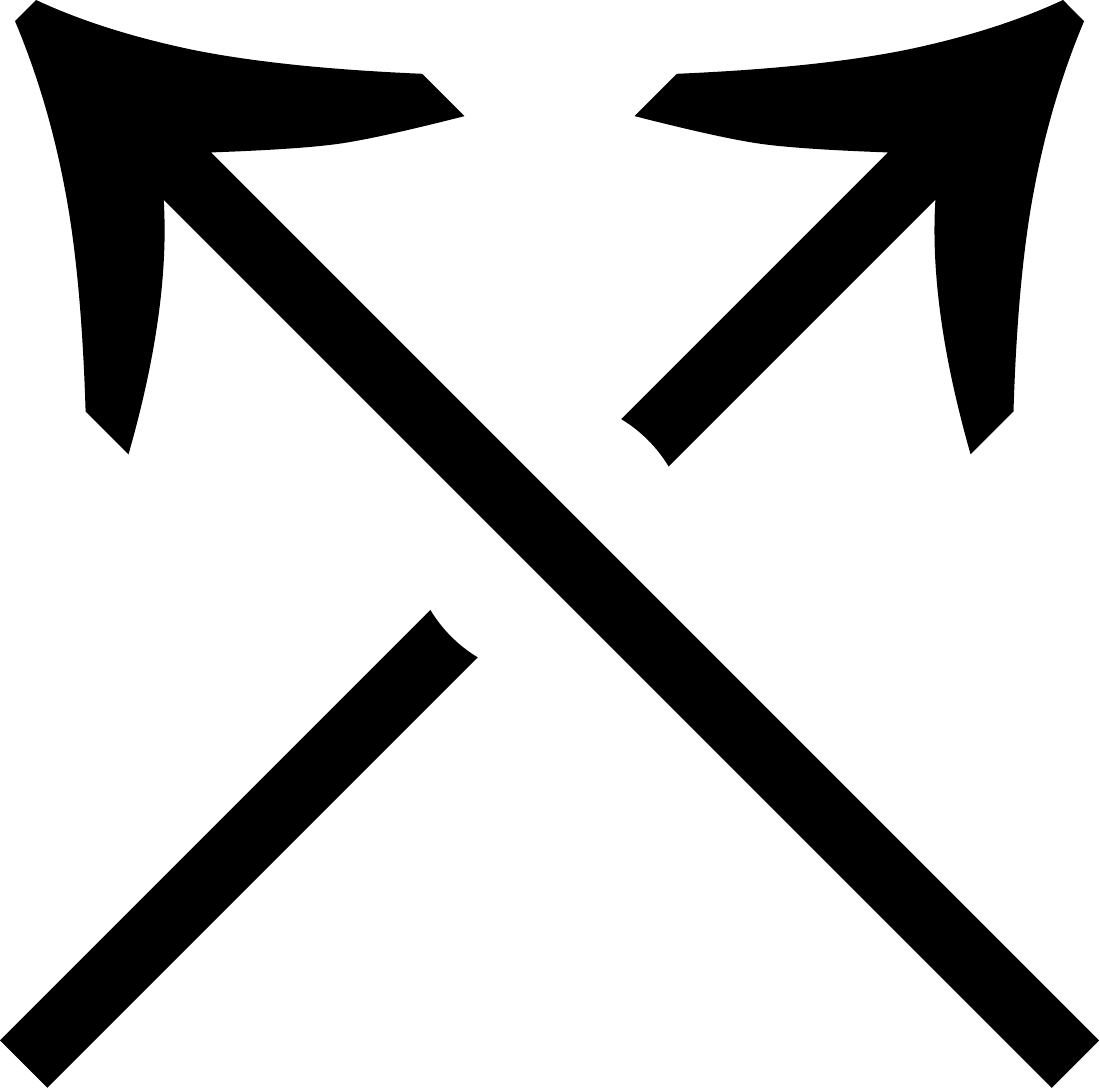}}%
  \end{picture}%
\endgroup%
}\vspace{-0.5em}
\begin{center}Right-handed\hspace{0.37\columnwidth}Left-handed\end{center}
\vspace{-1em}}
{The right-handed and left-handed crossings.
The crossings get their name from the fact that pointing the thumb of your right
hand along one of the strands in the right-handed crossing means your curled
fingers will be pointing in the direction of the other strand.
An analogous statement holds for left-handed crossings.
When computing values such as the writhe of a knot diagram, right-handed
crossings are considered positive and left-handed crossings are considered
negative.}

As described so far, this notation only tells us which strands cross which.
What it does not tell us is the handedness of each crossing (see \figCrossings).
In other words, the shape of the knot diagram can be reconstructed, but every
crossing will effectively be blurred out, as it will not be clear which of the
two strands in the crossing is the upper strand and which is the lower strand.
To account for this, we declare that a crossing is \textit{positive} if out of
the two values that make up a crossing, the odd one corresponds to the upper
strand of the crossing.
If a crossing is not positive, it is \textit{negative} and we indicate this by
negating the even value in each negative crossing.
For example, if a crossing is marked (17, 34) and the strand labeled 17 passes
above the strand labeled 34, we leave the crossing as is.
On the other hand, if the upper strand is marked 34, we denote the crossing by
the pair (17, -34).

If we were to flip over a knot diagram and look at it from the back, all of the
values in the MD code would change sign.
To account for this, when necessary we negate all of the values in the MD code
to make the leading term positive.
As a result, every knot diagram with $n$ crossings can be represented by a
signed permutation of the numbers from 1 to $n$.

\papersec{Alternating Knots}

A subset of the knots we are trying to tabulate are called
\textit{alternating knots}.
By determining which alternating knots should be included in our tabulation, we
can simplify the task of determining the remainder of the list.
For this reason, we start by determining which alternating knots should be
included in our reconstruction of the Rolfsen table.

When we move along a knot diagram, labeling its edges to determine its MD code,
we go over some strands and under others.
If we always alternate between going over and under the strands we cross, we say
that the knot diagram is \textit{alternating}.
We note that any minimal knot diagram that is equivalent to an alternating knot
diagram will be alternating (see \cite{new}).
Thus, it makes sense to refer to \textit{alternating knots} as this property is
independent of our choice of knot diagram, as long as it is minimal.
Any knot that is not an alternating knot is a \textit{non-alternating knot}.
We note that the MD code of an alternating knot will consist entirely of
positive entries.

To generate all of the knot diagrams that might be included in our tabulation,
we do not need to generate every possible knot diagram that there is.
It suffices to first determine the list of alternating knots in our
reconstruction of the Rolfsen table.
Afterwards, we will construct the non-alternating knots from our finalized list
of alternating knots.

We know that not all permutations of 10 values result in valid knots.
Thus, some of these permutations must be eliminated from consideration.
There are several criteria which we can use to determine which alternating knots
should be included in our list.
We will first define these criteria, then explain how to implement a test that
verifies that they are satisfied.
The alternating knot corresponding to a given permutation is included in our
tabulation if and only if it meets all of the following criteria.

\begin{enumerate}
\item A knot diagram can produce different permutations depending on where one
starts numbering and in which direction they proceed.
There are $4n$ ways to choose both a starting point and a direction.
A permutation is \textit{minimal} if it is lexicographically smaller than or
equal to all of the other $4n-1$ possible permutations of the corresponding knot
diagram.
To satisfy this criterion, a permutation must be minimal.

\item The resulting knot diagram must be prime.
Since a composite knot diagram can be split into two knot factors, we know that
as we label the knot, we will have to go through all of the crossings of one of
the knot factors before we move on to the other.
This means that the values from 1 to $2n$ will be split into two consecutive
subsequences since the values in each subsequence will be the labels of the
crossings of one of the knot factors.
In a permutation, this would be expressed as a set of $k$ pairs, all $2k$ of
whose values form a consecutive subsequence.
Thus, such a set must not exist in the ED code for the knot diagram to be prime
(see \figComposite).
This also handily eliminates knot diagrams that contain a kink and could be
simplified with the first Reidemeister move.
The third Reidemeister move and the simplifying direction of the second
Reidemeister move cannot occur in alternating knots (see \figMoves).

\item The permutation must encode a diagram which is realizable.
This means that there must be a way to draw the knot diagram in the plane
without adding any intersections beyond the ones encoded in the permutation.
The simplest permutation that fails this test is (2, 4, 1, 5, 3).
It is physically impossible to draw a knot diagram on a plane that would have a
non-realizable permutation as its MD code.

\item The knot diagram must be minimal with respect to flyping.
This means that the knot diagram's minimal permutation must be lexicographically
minimal over all of the permutations of knot diagrams that can be obtained from
our original diagram by applying a sequence of flypes (see \figMoves).
\end{enumerate}

The first two conditions can be used to avoid testing all $n!$ possible
alternating MD codes.
If we arranged the $n!$ permutations lexicographically and went along checking
each one, it will frequently be possible to skip checking up to $k!$
permutations at a time, where $0\leq k<n$.

Skipping these permutations is made possible as the first two criteria can be
checked directly from the permutation.
If a permutation that fails one of the first two criteria contains a string of
values that make the permutation fail this criterion, then all permutations that
contain the same string of values will also fail this criterion and do not need
to be considered.
Thus, as soon as we find such a permutation, we increment the last digit of the
offending subsequence, thereby skipping $k!$ permutations ahead, where $k$ is
the number of digits left in the permutation after the subsequence.
In other words, given a leading subsequence that fails a criterion, we do not
try to continue it with terms that we know will fail the same test.

For the first criterion, a permutation is not minimal if we can find a pair of
values in our ED code which are numerically closer together than the leading
term and its odd-valued pair, the number 1.
So if there is an even value $x$ in the ED code which is closer to its paired
odd value than the first number in the ED code is to 1, then all permutations
with $x$ in the same position will not be minimal.
This is because starting the enumeration of the knot diagram's crossings at the
one previously labeled $x$ would result in a smaller first element in the MD
code, which is not allowed as the knot diagrams in the table must be minimal.

\pdfmsgex{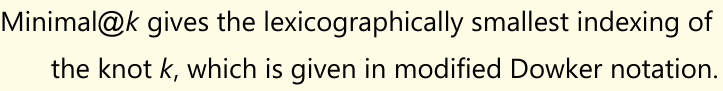}

Similarly, if a knot diagram is composite, this is represented by several
consecutive terms in the permutation, so all permutations obtained by
rearranging the values that come after this sequence would also fail this test.
Checking if a knot is prime was described above, we need to find a subsequence
of the values from 1 to $2n$ such that the all of the values' pairs are just a
reordering of the same subsequence.
For example, if an MD code starts with (2, 3, 1, \dots) it will not be prime.
This is because the subsequence (1, 2, 3, 4, 5, 6) has pairs (4, 5, 6, 1, 2, 3).
Since the second list is just a rearrangement of the first, any knot that starts
with this sequence will contain a knot factor, a trefoil in this case, and will
not be prime.

In other words, if a knot is composite, its knot factors will contain
consecutive strand labels.
Since two knot factors can never cross, the crossings of a knot factor contain
two elements of the subsequence that our permutation starts with.
Thus, if the pairs of the values in a subsequence do not contain any values
outside of that subsequence, the knot diagram is not prime.

The third condition is checked with the help of a modified graph planarity
algorithm.
If a 4-valent graph is constructed out of a knot diagram by replacing each
crossing with a vertex and each edge of the knot with an edge in the graph, then
typical planarity tests would frequently give false positives.
There are 4 edges emanating from a crossing, but there are only 2 ways of
arranging them in a valid manner in a knot diagram, but there are 6 ways of
arranging 4 edges around a vertex.
The reason for this is that a strand is not allowed to exit a crossing via an
edge that is adjacent to its incoming edge.
Strands must go directly across a crossing which means that the incoming and
outgoing edges of a strand must be aligned opposite from each other.

We have not yet imposed any restrictions that would tell a graph planarity
algorithm that such cases should not be considered.
Permutations that fail this test do not form a planar knot diagram, but the
graph that is created by making the same connections between vertices is
planar.
We can check that the graph formed by the permutation (2, 4, 1, 5, 3) is planar,
yet such a knot diagram is impossible to draw, and is thus not realizable.

\paperfig{Graph}{\begin{center}\def\svgwidth{0.8\columnwidth}
{\centering
\begingroup%
  \makeatletter%
  \providecommand\color[2][]{%
    \errmessage{(Inkscape) Color is used for the text in Inkscape, but the package 'color.sty' is not loaded}%
    \renewcommand\color[2][]{}%
  }%
  \providecommand\transparent[1]{%
    \errmessage{(Inkscape) Transparency is used (non-zero) for the text in Inkscape, but the package 'transparent.sty' is not loaded}%
    \renewcommand\transparent[1]{}%
  }%
  \providecommand\rotatebox[2]{#2}%
  \ifx\svgwidth\undefined%
    \setlength{\unitlength}{456.09172459bp}%
    \ifx\svgscale\undefined%
      \relax%
    \else%
      \setlength{\unitlength}{\unitlength * \real{\svgscale}}%
    \fi%
  \else%
    \setlength{\unitlength}{\svgwidth}%
  \fi%
  \global\let\svgwidth\undefined%
  \global\let\svgscale\undefined%
  \makeatother%
  \begin{picture}(1,0.36495715)%
    \put(0,0){\includegraphics[width=\unitlength,page=1]{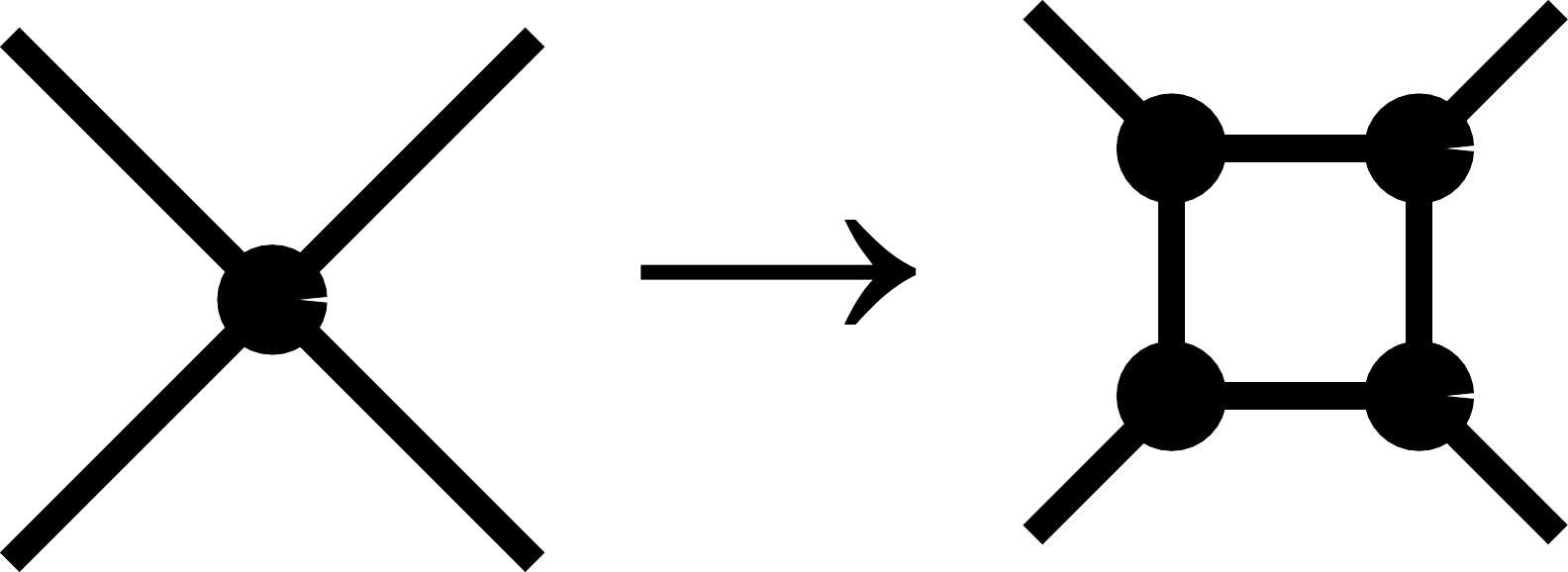}}%
  \end{picture}%
\endgroup%
}\end{center}
\vspace{-1em}}
{The transformation applied to a knot diagram's graph to determine if the knot
diagram is planar.
Each vertex in the 4-valent graph is replaced with 4 vertices connected to each
other and to the original edges in a square.
This makes the graph 3-valent and also serves as a proper indicator of the
planarity of the knot diagram's graph.
The reason for this is that a graph should not be accepted as planar if there is
a crossing where the two strands in the crossing enter and leave the crossing
through adjacent edges.
If this were to happen, the new graph would stop being planar as the square in
the centre would become a non-planar bowtie.}

To solve this problem, it is sufficient to replace each vertex with four
vertices in a square to construct the modified graph of the knot diagram (see
\figGraph).
This preserves the planarity of the two allowable configurations but bars the
other four, as the square would be transformed into a non-planar bowtie shape.
Thus, it suffices to use a regular graph planarity algorithm to check whether
the modified graph of the knot diagram is planar.
If it is not, the permutation fails the third criterion.

\pdfmsgex{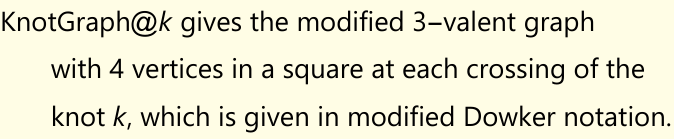}

Finally, the fourth condition is checked by using a graph searching algorithm to
find all knot diagrams that can be obtained from a given knot diagram with a
sequence of flypes (see \figMoves).
From these, we keep only the lexicographically minimal knot diagram.
A flype is represented in a permutation as a pair and two disjoint subsequences
of 1 to $2n$.
The two subsequences are the strands that are the part of the knot that gets
flipped and the pair is the crossing that gets moved to the other side of the
knot during the flype.
Note that these subsequences must satisfy several conditions.
The first is that all of the numbers in these subsequences must be pairs of
other elements of the subsequences.
This is much like testing a knot for primality.
The reason for this condition is that the part of the knot that is flyped must
be like a two strand knot factor, in the sense that it must not cross any part
of the knot outside of itself.
Additionally, these two subsequences, depending on which way the strands run,
must either start or end with the values adjacent to those of the earlier pair.
Since that pair contains one odd value and one even value, we can save time by
ignoring the cases where this would be violated.

\pdfmsgex{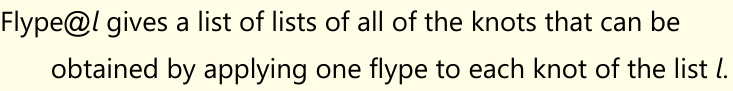}

Without the flyping condition, we have generated a list of all candidate knot
diagrams.
This is a list of all diagrams that encode alternating knots and all appear
distinct at first glance.
We will use this list later when we will want to examine all valid knot diagrams
as opposed to those that are necessarily distinct.

\pdfmsgex{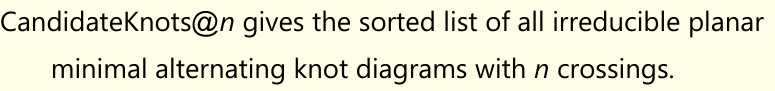}

At this point, it is reasonable to doubt that the four conditions have
completely narrowed down the list of alternating knots, leaving no equivalent
knot diagrams.
The reason for why this is sufficient is because any two minimal equivalent
alternating knot diagrams are related by a sequence of flypes (see
\cite{flype}).
This means that our list of alternating knot diagrams contains no duplicates and
is complete.
There are 197 alternating knots with 10 crossings or fewer.

\pdfmsgex{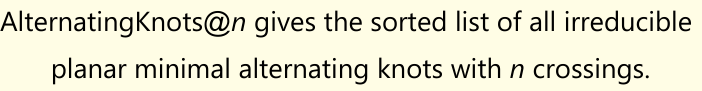}

\papersec{Non-Alternating Knots}

After generating all of the alternating knots, we have a smaller set of knot
diagrams to test to see if they should be included in our tabulation.
There are $2^{n-1}n!$ signed permutations with a positive leading term.
Of these, we only need to consider those that have the shape of alternating
knots in our list.
What this means is that we only have to consider the $2^{n-1}$ knot diagrams
obtained for each alternating knot in our list by flipping the signs of the
elements of the alternating knot's MD code in every possible way.
There would be $2^n$ ways to do this but the leading term must stay positive,
leaving us with $2^{n-1}$ ways.

For the case where $n=10$, we have 123 alternating knots.
This means that we will generate 62976 non-alternating knots.
Although this may seem like a large number of knots, most of these can be
examined and discarded immediately.

Almost all of these knot diagrams are discarded as they can be reduced with the
second Reidemeister move.
Many of those that remain can be reduced with a (2, 1)--pass or a (3, 2)--pass
(see \figMoves).
We do not consider the (1, 0)--pass as it is the first Reidemeister move.
These pass moves can be found in most reducible knot diagrams.

\pdfmsgex{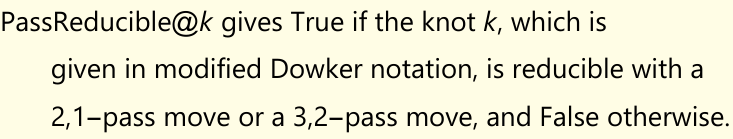}

By considering all knot diagrams that cannot be simplified with either the
second Reidemeister move or a pass move, there are very few reducible knot
diagrams that have not yet been eliminated.
The reason for this is that a reducibility test checking only for these moves
will occasionally give false negatives.
This will be dealt with at a later stage.
We will call knot diagrams that can be reduced with the second Reidemeister move
or a pass move as \textit{immediately reducible}.
Note that not all reducible knots are immediately reducible.
At this point, we have 1176 non-alternating knot diagrams that are not
immediately reducible left to consider.

\pdfmsgex{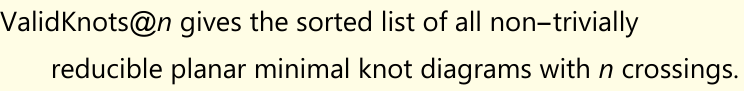}

\papersec{Finding Equivalent Diagrams}

Canonically, alternating knots precede non-alternating knots in the Rolfsen
table.
We maintain the same pattern, ordering the knots within each of the two sets
lexicographically.
Additionally, knots with a lower crossing number always precede those with a
higher crossing number.

\pdfmsgex{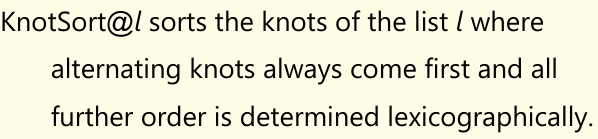}

We need to determine which knot diagrams are equivalent and find the
lexicographically smallest permutation for each knot.
We do this by examining all knot diagrams which cannot be transformed into
lexicographically smaller ones by applying a sequence of crossing
number-preserving moves.

For 10 crossings and fewer, the third Reidemeister move, the 2--pass, and the
flype (see \figMoves) are sufficient to eliminate all but 54 of the
non-alternating knot diagrams that are not immediately reducible.
Applying the third Reidemeister move is preferable to applying either of the
first two Reidemeister moves.
This is due to the fact that there are many ways of adding crossings to a knot
diagram but there are only a few ways to apply a move that preserves the knot
diagram's crossing number.

\pdfmsgex{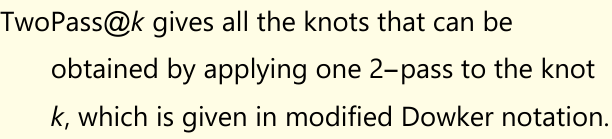}

To find equivalent knot diagrams, we implement a graph searching algorithm.
We first need to build our graph recursively by adding on subsequent layers of
knot diagrams.
We start with a graph $\Gamma_0$ consisting of the set of vertices $V_0$ and
edges $E_0$.
We define $V_0$ as the set of those 1176 knot diagrams that are not immediately
reducible.
For every natural number $i$, we define $V_i$ as the union of $V_{i-1}$ and the
set of knot diagrams that can be obtained by applying one crossing
number-preserving move to a knot diagram in $V_{i-1}$.
We do not need to define our edges recursively.
For any non-negative integer $i$, we define $E_i\subset V_i\times V_i$.
For any two knot diagrams $v_1$ and $v_2$ in $V_i$ we include the undirected
edge ($v_1$, $v_2$) in $E_i$ if we can apply a crossing number-preserving move
to $v_1$ and obtain $v_2$ (see \figMoves).
Lastly, $\Gamma_i$ is simply the set of vertices $V_i$ and set of edges $E_i$.

Since there are finitely many knot diagrams with a given number of crossings,
there must exist an integer $i$ such that $\Gamma_i$ is equivalent to
$\Gamma_{i+1}$, at which point the graph will cease to change.
We then take $\Gamma_i$ to be our graph.
Each connected component of the graph consists of a set of equivalent diagrams,
all representing the same knot.

At this point we return to the earlier concern that this graph might contain
some reducible knot diagrams.
Any reducible knot diagram that has not yet been removed is not immediately
reducible, which means that it cannot be reduced with a pass move or the second
Reidemeister move.
However, all such diagrams are equivalent to other diagrams which are
immediately reducible.
Thus, we check to see if any of the knot diagrams in a connected component are
immediately reducible.
If at least one is, we need to remove the entire component from the graph.
To do this, we create a graph $\Gamma'$ which is a subgraph of $\Gamma_i$ and
contains only the edges and vertices of the components which do not contain any
immediately reducible knot diagrams.
After this, $\Gamma'$ does not contain any reducible knot diagrams.

We now must generate our list of knots from the graph $\Gamma'$.
To do this, we take the lexicographically smallest knot diagram from each
connected component of $\Gamma'$.
As previously mentioned, we get 54 such knots.
Combined with the 197 alternating knots, this gives us 251 total knots.
However, just because we have applied a variety of moves to construct edges in
our graph does not mean that we are done.
It is possible that there are equivalent knot diagrams in the graph $\Gamma'$
between which we were unable to find a sequence of moves out of our set.
Thus, all we know is that there are no more than 251 knots with 10 crossings or
fewer.
Our lower bound is currently 200, as we know that our alternating knots are
distinct and that we have at least one non-alternating knot with each of 8, 9,
and 10 crossings.

The reason that we may have more knot diagrams than we should is because we are
restricting ourselves to using flypes, 2--passes, and the third Reidemeister
move to find equivalent knot diagrams.
Reidemeister's original theorem has the consequence that if two minimal knot
diagrams are equivalent and it is impossible to transform one into the other by
repeatedly applying the third Reidemeister move, the only option left available
to us is to first add crossings using one of the first two Reidemeister moves,
and proceed from there.
We are in a similar position because to show that two of our knot diagrams are
equivalent, we would have to turn them into more complicated equivalent knot
diagrams.
However, though there are relatively few ways to apply the third Reidemeister
move to a knot diagram, there are many ways of adding a kink to one.

Our next step is to figure out which pairs of knot diagrams in our list might be
equivalent.
Since there are so many ways of adding kinks, checking all of our knot diagrams
for equivalence this way would take a long time.
Thus, we would like some way to establish with certainty that some of our knot
diagrams cannot be equivalent to any others.
We do this by using invariants.

\papersec{Invariants}

All invariants are functions from the space of knot diagrams to an arbitrary
target space.
The property that an invariant must satisfy is that the images of two equivalent
knot diagrams must be equal under an invariant.
This allows us to show that two knot diagrams on which an invariant produces a
different value are knot diagrams of two different knots.

The invariant condition is easy to satisfy as we could choose the image of the
invariant to have unit magnitude.
For example, we could state that for any knot diagram, our carefully crafted
invariant produces a value of 0.
However, such an invariant is useless as we want to be able to show that some
knot diagrams are distinct.
Thus, we need invariants that produce the same value for non-equivalent knot
diagrams as rarely as possible.
The degree to which an invariant accomplishes this it is typically called its
strength, where weak invariants often produce the same values for different knot
diagrams and vice versa.
Often, stronger invariants require more time to compute, which is why it is
useful to have several invariants.
We use them in increasing order of strength so that the most difficult
computations only have to be done for a few knot diagrams, those between which
the weaker invariants were unable to distinguish.

The weakest invariant that we have available to us is crossing number.
We know that all of our knots are non-reducible and we also know that all of our
alternating knots are distinct.
This immediately reduces our task to simply calculating invariants on
non-alternating knots that all have the same crossing number.

After this initial step, we use two different invariants to complete the task.
Our initial simplifications would typically be considered too crude to be deemed
invariants except in the most technical of circumstances.
Our first invariant, the Jones polynomial, is fast and fairly strong and the
second, the number of colourings, is slow but even stronger.

\pdfmsgex{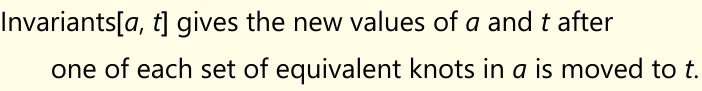}

\papersec{Planar Diagram Notation}

To find the Jones polynomial of a given knot diagram, the knot diagram must be
written in planar diagram notation as opposed to as an MD code.
To find this notation, it is necessary to determine the handedness of each of
the knot diagram's crossings (see \figCrossings).
All we know about a knot diagram's crossings from an MD code is which strand
passes above or below.
We do not know the handedness of each crossing.
As there are 2 types of crossings, there are $2^n$ possible ways to set the
handedness of a knot diagram's crossings.
Since the only knot diagrams being considered are realizable, it is known that
at least one of these $2^n$ crossing orientations will make the knot diagram
planar.

Since we are trying to compute polynomials of knot diagrams where $n\leq10$,
we have that $2^n\leq1024$, which is, computationally speaking, a small number.
For this reason, we can exhaustively iterate through the $2^n$ crossing
orientations until we find one that creates a planar knot diagram.

\papersvg{Ribbon}{ribbon}
{The transformation applied to the knot diagram's graph to determine whether or
not the given configuration of crossings makes the knot diagram planar.
We take each vertex in the 4-valent graph and replace it with 4 vertices
connected to each other and to the original edges in a diamond.
However, the connections to the adjacent vertices have been expanded to be a
pair of parallel strands.
This serves as a proper indicator of the planarity of the knot diagram's graph
with the given crossing configuration.
The reason for this is that a graph should not be accepted as planar if there is
a crossing where the two strands in the crossing enter and leave the crossing
through adjacent edges.
The new graph would stop being planar if this were to happen since the diamond
in the centre would become a non-planar bowtie.
Additionally, the graph should not be accepted as planar if the handedness of
the crossing is changed.
If this happens, the ribbons would twist and stop being planar.}

To test if the knot diagram with given crossing orientations is planar, we
apply the knot diagram planarity replacement from before, but replace each outer
edge with a ribbon, a pair of parallel edges, to only allow 1 of the 6 edge
configurations (see \figRibbon).

Earlier, the square replacement remained planar for either of the two ways of
arranging the edges of the crossing so that the strands enter and exit the
crossing through opposite edges (see \figGraph).
Now, we wish to exclude one of these two configurations.
We do this by arranging the strands into the configuration we desire and
changing the incoming edges into pairs of edges.
Now, whenever the strands do not exit and enter through opposite edges, the
graph will not be planar, just as before and for the same reasons.
More importantly, when the handedness of the crossing changes, a twist will be
added to two of the ribbons, making the graph non-planar.
Thus, all we need to do is search through all possible sets of crossing
orientations until we find one for which this modified graph is planar.

\paperfig{X}{
\begin{center}$k$\hspace{0.2\columnwidth}$j$\\
\def\svgwidth{0.2\columnwidth}
{\centering
\begingroup%
  \makeatletter%
  \providecommand\color[2][]{%
    \errmessage{(Inkscape) Color is used for the text in Inkscape, but the package 'color.sty' is not loaded}%
    \renewcommand\color[2][]{}%
  }%
  \providecommand\transparent[1]{%
    \errmessage{(Inkscape) Transparency is used (non-zero) for the text in Inkscape, but the package 'transparent.sty' is not loaded}%
    \renewcommand\transparent[1]{}%
  }%
  \providecommand\rotatebox[2]{#2}%
  \ifx\svgwidth\undefined%
    \setlength{\unitlength}{316.54126374bp}%
    \ifx\svgscale\undefined%
      \relax%
    \else%
      \setlength{\unitlength}{\unitlength * \real{\svgscale}}%
    \fi%
  \else%
    \setlength{\unitlength}{\svgwidth}%
  \fi%
  \global\let\svgwidth\undefined%
  \global\let\svgscale\undefined%
  \makeatother%
  \begin{picture}(1,0.98976115)%
    \put(0,0){\includegraphics[width=\unitlength,page=1]{diagrams/positive.pdf}}%
  \end{picture}%
\endgroup%
}\\
$l$\hspace{0.2\columnwidth}$i$\\
$X_{i,j,k,l}$\end{center}\vspace{-1em}}
{A right-handed crossing labeled in planar diagram notation.
The lower incoming edge is labeled $i$ and then the remaining three are labeled
$j$, $k$, and $l$, proceeding counterclockwise from $i$.
The crossing is labeled as $X_{i,j,k,l}$.}

Every crossing is represented in planar diagram notation as $X_{i,j,k,l}$ (see
\figX).
Here, $i$ is the index given to the lower incoming edge and then $j$, $k$, and
$l$ proceed counterclockwise.
The knot diagram is then written as the product of its crossings in planar
diagram notation.
For example, the left-handed trefoil (see \figLabeled) is written as
$X_{1,4,2,5}X_{3,6,4,1}X_{5,2,6,3}$.

\pdfmsgex{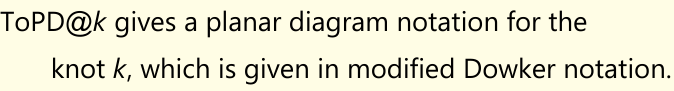}

Once we can transform knot diagrams into planar diagram notation, we can compute
their Jones polynomials.

\papersec{Jones Polynomial}

The Jones polynomial of a knot diagram is computed from the product of the knot
diagram's crossings.\\

\paperfig{Smoothings}{
\begin{center}Smoothings for a Right-Handed Crossing\end{center}
\def\svgwidth{0.25\columnwidth}
{\centering
\begingroup%
  \makeatletter%
  \providecommand\color[2][]{%
    \errmessage{(Inkscape) Color is used for the text in Inkscape, but the package 'color.sty' is not loaded}%
    \renewcommand\color[2][]{}%
  }%
  \providecommand\transparent[1]{%
    \errmessage{(Inkscape) Transparency is used (non-zero) for the text in Inkscape, but the package 'transparent.sty' is not loaded}%
    \renewcommand\transparent[1]{}%
  }%
  \providecommand\rotatebox[2]{#2}%
  \ifx\svgwidth\undefined%
    \setlength{\unitlength}{333.93346956bp}%
    \ifx\svgscale\undefined%
      \relax%
    \else%
      \setlength{\unitlength}{\unitlength * \real{\svgscale}}%
    \fi%
  \else%
    \setlength{\unitlength}{\svgwidth}%
  \fi%
  \global\let\svgwidth\undefined%
  \global\let\svgscale\undefined%
  \makeatother%
  \begin{picture}(1,1.00148955)%
    \put(0,0){\includegraphics[width=\unitlength,page=1]{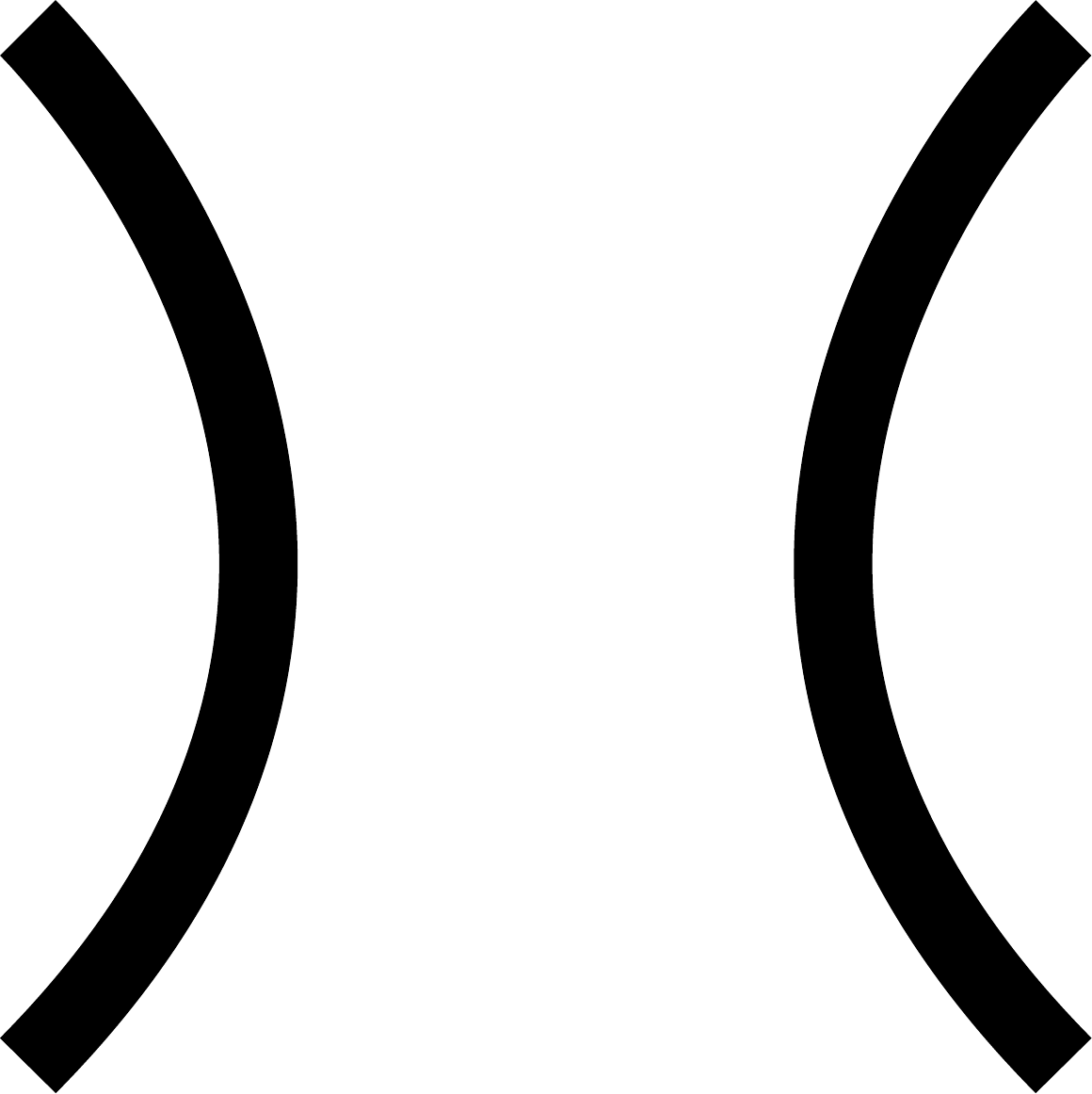}}%
  \end{picture}%
\endgroup%
}
\hfill
\def\svgwidth{0.25\columnwidth}
{\centering
\begingroup%
  \makeatletter%
  \providecommand\color[2][]{%
    \errmessage{(Inkscape) Color is used for the text in Inkscape, but the package 'color.sty' is not loaded}%
    \renewcommand\color[2][]{}%
  }%
  \providecommand\transparent[1]{%
    \errmessage{(Inkscape) Transparency is used (non-zero) for the text in Inkscape, but the package 'transparent.sty' is not loaded}%
    \renewcommand\transparent[1]{}%
  }%
  \providecommand\rotatebox[2]{#2}%
  \ifx\svgwidth\undefined%
    \setlength{\unitlength}{334.50384896bp}%
    \ifx\svgscale\undefined%
      \relax%
    \else%
      \setlength{\unitlength}{\unitlength * \real{\svgscale}}%
    \fi%
  \else%
    \setlength{\unitlength}{\svgwidth}%
  \fi%
  \global\let\svgwidth\undefined%
  \global\let\svgscale\undefined%
  \makeatother%
  \begin{picture}(1,1.00000597)%
    \put(0,0){\includegraphics[width=\unitlength,page=1]{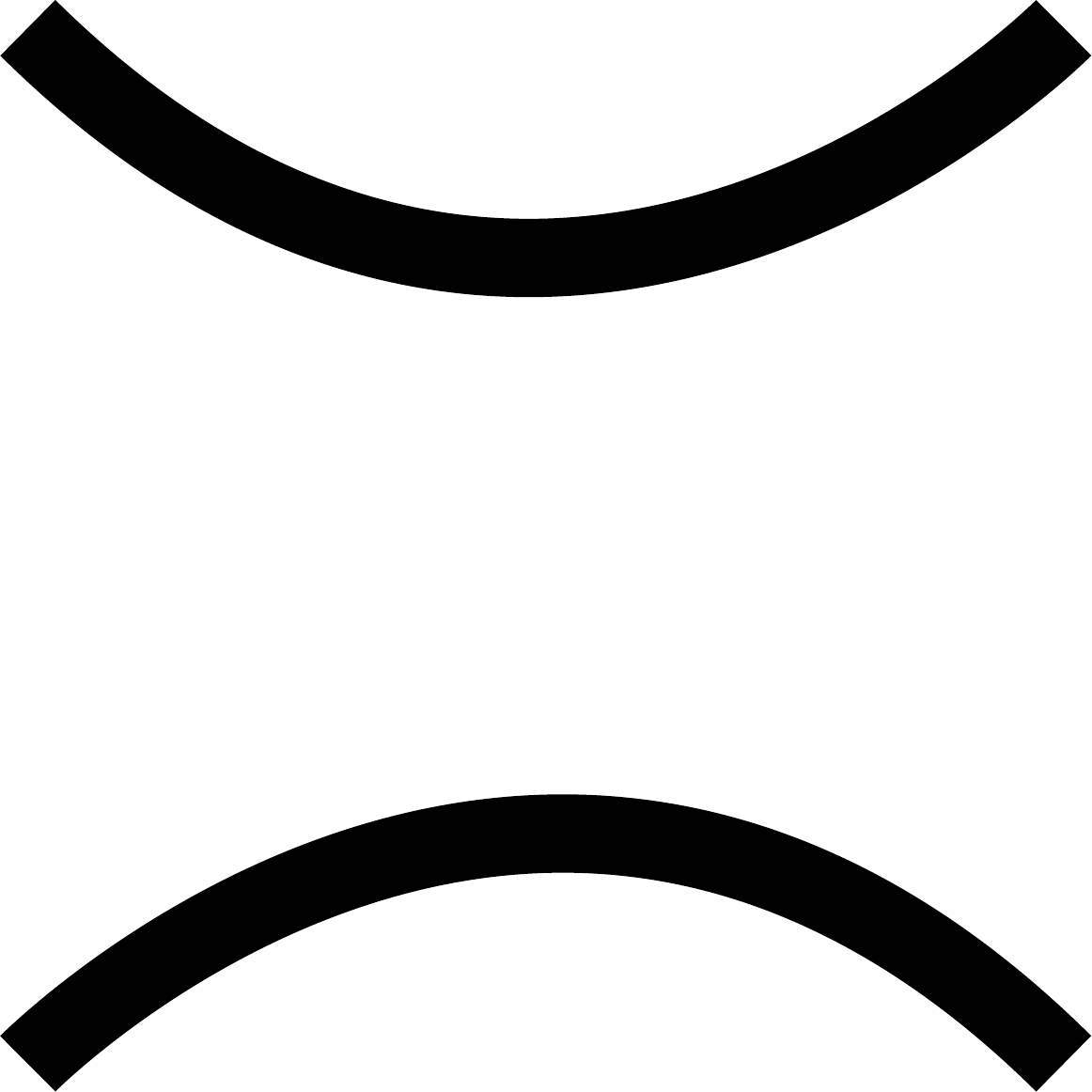}}%
  \end{picture}%
\endgroup%
}\\

\noindent0-smoothing\hfill 1-smoothing}
{The 0 and 1-smoothings of a right-handed crossing.
The smoothings are comprised of two strands with no directionality.
If every crossing in a knot diagram is replaced by a smoothing, the result is an
unlink as the knot diagram will be devoid of any crossings or ends.
The 0-crossing is formed by connecting each of the two ends of the lower strand
of the crossing to the ends of the upper strand that are next to them in the
counterclockwise direction.
For the 1-smoothing, the direction is clockwise.
The 0-smoothing for a right-handed crossing is identical to the 1-smoothing for
a left-handed crossing and vice versa.}

Every crossing can be \textit{smoothed} in two distinct ways (see
\figSmoothings).
By smoothing a crossing in a particular manner, the polynomial of that smoothing
is multiplied by a coefficient of either $A$ or $B$ for the 0 and 1-smoothings,
respectively.
Since each smoothing is actually a coefficient multiplied by the two
non-intersecting strands of the smoothing, a strand stitching operation is
applied to turn a product of $n$ smoothings into an unlink of several
components.

A strand stitching operation satisfies the property that the product of two
strands that share an endpoint, such as the strand from $p$ to $q$, ($p$, $q$),
and the strand ($q$, $r$), will be equal to one strand running between their
non-common endpoints, ($p$, $r$) in this case.

The final result will always be the product of several strands that are closed
loops of the form ($p$, $p$).
Each of these components of the link is given a coefficient of $d$ and thus the
result is a polynomial in $A$, $B$, and $d$.
What we have defined so far is called the Kauffman bracket of a knot diagram
$X$, and it is denoted $\langle X\rangle$.
We note that $\langle\bigcirc\rangle=d$ and $\langle$\o$\rangle=1$, where
$\bigcirc$ and \o~represent the unknot and the empty knot, respectively.
Using this notation, a formula for the smoothings of a crossing can be written.

\papereq{BracketPlus}{
\left\langle\begin{matrix}\def\svgwidth{2em}
{\centering
\begingroup%
  \makeatletter%
  \providecommand\color[2][]{%
    \errmessage{(Inkscape) Color is used for the text in Inkscape, but the package 'color.sty' is not loaded}%
    \renewcommand\color[2][]{}%
  }%
  \providecommand\transparent[1]{%
    \errmessage{(Inkscape) Transparency is used (non-zero) for the text in Inkscape, but the package 'transparent.sty' is not loaded}%
    \renewcommand\transparent[1]{}%
  }%
  \providecommand\rotatebox[2]{#2}%
  \ifx\svgwidth\undefined%
    \setlength{\unitlength}{316.54126374bp}%
    \ifx\svgscale\undefined%
      \relax%
    \else%
      \setlength{\unitlength}{\unitlength * \real{\svgscale}}%
    \fi%
  \else%
    \setlength{\unitlength}{\svgwidth}%
  \fi%
  \global\let\svgwidth\undefined%
  \global\let\svgscale\undefined%
  \makeatother%
  \begin{picture}(1,0.98976115)%
    \put(0,0){\includegraphics[width=\unitlength,page=1]{diagrams/positive.pdf}}%
  \end{picture}%
\endgroup%
}\end{matrix}\right\rangle
=A\left\langle\begin{matrix}\def\svgwidth{2em}
{\centering
\begingroup%
  \makeatletter%
  \providecommand\color[2][]{%
    \errmessage{(Inkscape) Color is used for the text in Inkscape, but the package 'color.sty' is not loaded}%
    \renewcommand\color[2][]{}%
  }%
  \providecommand\transparent[1]{%
    \errmessage{(Inkscape) Transparency is used (non-zero) for the text in Inkscape, but the package 'transparent.sty' is not loaded}%
    \renewcommand\transparent[1]{}%
  }%
  \providecommand\rotatebox[2]{#2}%
  \ifx\svgwidth\undefined%
    \setlength{\unitlength}{333.93346956bp}%
    \ifx\svgscale\undefined%
      \relax%
    \else%
      \setlength{\unitlength}{\unitlength * \real{\svgscale}}%
    \fi%
  \else%
    \setlength{\unitlength}{\svgwidth}%
  \fi%
  \global\let\svgwidth\undefined%
  \global\let\svgscale\undefined%
  \makeatother%
  \begin{picture}(1,1.00148955)%
    \put(0,0){\includegraphics[width=\unitlength,page=1]{diagrams/zero.pdf}}%
  \end{picture}%
\endgroup%
}\end{matrix}\right\rangle
+B\left\langle\begin{matrix}\def\svgwidth{2em}
{\centering
\begingroup%
  \makeatletter%
  \providecommand\color[2][]{%
    \errmessage{(Inkscape) Color is used for the text in Inkscape, but the package 'color.sty' is not loaded}%
    \renewcommand\color[2][]{}%
  }%
  \providecommand\transparent[1]{%
    \errmessage{(Inkscape) Transparency is used (non-zero) for the text in Inkscape, but the package 'transparent.sty' is not loaded}%
    \renewcommand\transparent[1]{}%
  }%
  \providecommand\rotatebox[2]{#2}%
  \ifx\svgwidth\undefined%
    \setlength{\unitlength}{334.50384896bp}%
    \ifx\svgscale\undefined%
      \relax%
    \else%
      \setlength{\unitlength}{\unitlength * \real{\svgscale}}%
    \fi%
  \else%
    \setlength{\unitlength}{\svgwidth}%
  \fi%
  \global\let\svgwidth\undefined%
  \global\let\svgscale\undefined%
  \makeatother%
  \begin{picture}(1,1.00000597)%
    \put(0,0){\includegraphics[width=\unitlength,page=1]{diagrams/one.pdf}}%
  \end{picture}%
\endgroup%
}\end{matrix}\right\rangle}{}
\papereq{BracketMinus}{
\left\langle\begin{matrix}\def\svgwidth{2em}
{\centering
\begingroup%
  \makeatletter%
  \providecommand\color[2][]{%
    \errmessage{(Inkscape) Color is used for the text in Inkscape, but the package 'color.sty' is not loaded}%
    \renewcommand\color[2][]{}%
  }%
  \providecommand\transparent[1]{%
    \errmessage{(Inkscape) Transparency is used (non-zero) for the text in Inkscape, but the package 'transparent.sty' is not loaded}%
    \renewcommand\transparent[1]{}%
  }%
  \providecommand\rotatebox[2]{#2}%
  \ifx\svgwidth\undefined%
    \setlength{\unitlength}{316.54130073bp}%
    \ifx\svgscale\undefined%
      \relax%
    \else%
      \setlength{\unitlength}{\unitlength * \real{\svgscale}}%
    \fi%
  \else%
    \setlength{\unitlength}{\svgwidth}%
  \fi%
  \global\let\svgwidth\undefined%
  \global\let\svgscale\undefined%
  \makeatother%
  \begin{picture}(1,0.9897609)%
    \put(0,0){\includegraphics[width=\unitlength,page=1]{diagrams/negative.pdf}}%
  \end{picture}%
\endgroup%
}\end{matrix}\right\rangle
=A\left\langle\begin{matrix}\def\svgwidth{2em}
{\centering
\begingroup%
  \makeatletter%
  \providecommand\color[2][]{%
    \errmessage{(Inkscape) Color is used for the text in Inkscape, but the package 'color.sty' is not loaded}%
    \renewcommand\color[2][]{}%
  }%
  \providecommand\transparent[1]{%
    \errmessage{(Inkscape) Transparency is used (non-zero) for the text in Inkscape, but the package 'transparent.sty' is not loaded}%
    \renewcommand\transparent[1]{}%
  }%
  \providecommand\rotatebox[2]{#2}%
  \ifx\svgwidth\undefined%
    \setlength{\unitlength}{334.50384896bp}%
    \ifx\svgscale\undefined%
      \relax%
    \else%
      \setlength{\unitlength}{\unitlength * \real{\svgscale}}%
    \fi%
  \else%
    \setlength{\unitlength}{\svgwidth}%
  \fi%
  \global\let\svgwidth\undefined%
  \global\let\svgscale\undefined%
  \makeatother%
  \begin{picture}(1,1.00000597)%
    \put(0,0){\includegraphics[width=\unitlength,page=1]{diagrams/one.pdf}}%
  \end{picture}%
\endgroup%
}\end{matrix}\right\rangle
+B\left\langle\begin{matrix}\def\svgwidth{2em}
{\centering
\begingroup%
  \makeatletter%
  \providecommand\color[2][]{%
    \errmessage{(Inkscape) Color is used for the text in Inkscape, but the package 'color.sty' is not loaded}%
    \renewcommand\color[2][]{}%
  }%
  \providecommand\transparent[1]{%
    \errmessage{(Inkscape) Transparency is used (non-zero) for the text in Inkscape, but the package 'transparent.sty' is not loaded}%
    \renewcommand\transparent[1]{}%
  }%
  \providecommand\rotatebox[2]{#2}%
  \ifx\svgwidth\undefined%
    \setlength{\unitlength}{333.93346956bp}%
    \ifx\svgscale\undefined%
      \relax%
    \else%
      \setlength{\unitlength}{\unitlength * \real{\svgscale}}%
    \fi%
  \else%
    \setlength{\unitlength}{\svgwidth}%
  \fi%
  \global\let\svgwidth\undefined%
  \global\let\svgscale\undefined%
  \makeatother%
  \begin{picture}(1,1.00148955)%
    \put(0,0){\includegraphics[width=\unitlength,page=1]{diagrams/zero.pdf}}%
  \end{picture}%
\endgroup%
}\end{matrix}\right\rangle}{}

A given smoothing is a 0-smoothing if the incoming end of the lower strand is
connected to the next end going counterclockwise around the crossing, in other
words, the nearest end on its right.
If it is connected to the end on its left, the resulting smoothing is a
1-smoothing.

Thus, $\langle\svgl{right}\rangle$, the Kaufmann bracket of the right-handed
trefoil can be evaluated.
Note that this trefoil is right-handed so we will only need \eqBracketPlus.
There are two ways to smooth each crossing so there are eight ways to smoooth
the three crossings altogether.
Two of these ways are applying three 0-smoothings and three 1-smoothings.
The other six cases are not all distinct, since there are three identical ways
of applying either one or two 0-smoothings.
Thus, each of the cases in these sets have the same bracket value, which is how
the bracket of the trefoil is expanded.

\papereq{TrefoilOne}{\fontsize{9pt}{1em}\selectfont
\langle\svgl{right}\rangle
=A^3\langle\svgl{triple}\rangle
+3A^2B\langle\svgl{double}\rangle
+3AB^2\langle\svgl{single}\rangle
+B^3\langle\svgl{nil}\rangle}{}

Unsurprisingly, \eqTrefoilOne looks a lot like an application of the binomial
theorem.
However, that is only because the trefoil is rotationally symmetric.
In the general case, the bracket will not have as many like terms.

By counting the number of components in each unlink, the remaining
brackets are evaluated with the corresponding power of $d$.

\papereq{TrefoilTwo}{\fontsize{11pt}{1em}\selectfont
\langle\svgl{right}\rangle
=A^3d^2+3A^2Bd+3AB^2d^2+B^3d^3}{}

To make the Kaufmann bracket invariant over the second Reidemeister move, we
need to set $d+A^2+B^2=0$ and $AB=1$.
We find that these relations make the Kaufmann bracket invariant over the third
move as well.
This means that to show that this is actually an invariant, it suffices to show
that the Kaufamnn bracket is invariant over the first Reidemeister move.
We find that in its current form, this is not the case.

\pdfmsg{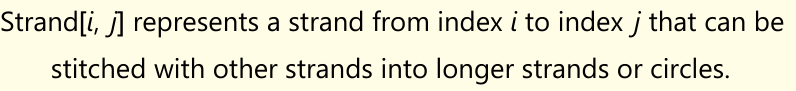}

To make the Kaufmann bracket invariant over addition and removal of kinks, the
whole polynomial must be multiplied by a coefficient of $(\text-A)^{\text-3w}$
where $w$ is the writhe of the knot diagram, which is the difference between the
number of right-handed and left-handed crossings in the knot diagram.

\pdfmsgex{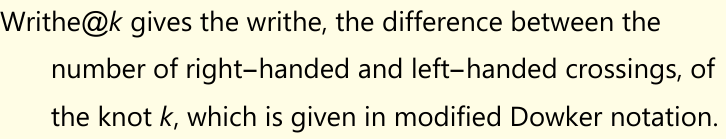}

The resulting polynomial will have a factor of $d$ in every component as every
unlink that we can get by smoothing the knot diagram has at least one component.
Thus, the polynomial is normalized by dividing it by $d$.
Lastly, what we have obtained will always be a polynomial in $A^4$ so the rule
$A=q^{\text-1/4}$ is applied to make the result a Laurent polynomial in $q$
(see \cite{new}).

For the trefoil, these substitutions allow us to transform our equation into a
simpler form.
We get that the writhe, $w$, is equal to 3.
This means that the Jones polynomial for the right-handed trefoil needs to be
multiplied by $\text-A^{\text-9}$.
Applying $d=\text-A^2-B^2$ and $B=A^{\text-1}$, we get the Jones polynomial of
the trefoil.

\papereq{TrefoilThree}{J(\svgl{right})
=\text-A^{\text-16}+A^{\text-12}+A^{\text-4}}{}

Since the Jones polynomial of the mirror image of a knot diagram is the Jones
polynomial of the original knot diagram with $q$ replaced by $q^{\text-1}$, the
minimal of these two polynomials is taken as the value of the invariant for that
knot diagram.

Applying the $q$ substitution will yield the final version of the Jones
polynomial for the right-handed trefoil.
However, the left-handed trefoil has a smaller Jones polynomial by degree so we
state that the Jones polynomial for the trefoil is the Jones polynomial for the
left-handed trefoil.

\papereq{TrefoilJones}{J(\svgl{left})
=\text-q^{\text-4}+q^{\text-3}+q^{\text-1}}{}

\pdfmsgex{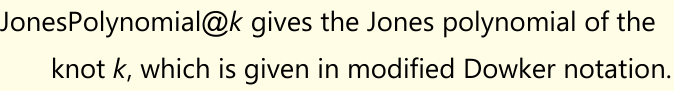}

We find that the Jones polynomial shows that every knot diagram out of our 251
with 9 crossings or fewer is distinct.
This is because all of the non-alternating knots with fewer than 10 crossings
have distinct values for their Jones polynomial.
Among the knots with 10 crossings, we find two pairs of knot diagrams with the
same Jones polynomial.
As we have dramatically reduced the list of diagrams we are unsure about, we can
now apply our more powerful invariant at little cost.

Note that this raises our lower bound to 249.
This is because there can be at most two extra knots in our list as the Jones
polynomial only found two pairs of knots whose Jones polynomial was not
distinct.

\papersec{Knot Colourings}

We have two pairs of knots in our list that could be equivalent.
We need to determine if the knot diagrams in either pair are distinct.
As we have very few knot diagrams to analyze, we can spend some additional time
computing a more complicated invariant, in exchange for it being able to
distinguish between our knot diagrams.
This invariant is the number of \textit{colourings} of the knot diagram with
elements of the permutation group $S_m$, for some $m$.
Such a colouring is an assignment of permutations of $m$ elements to edges of
the knot diagram such that these permutations satisfy a particular set of
conditions.
The number of such colourings is invariant over the three Reidemeister moves,
making it invariant over all equivalent knot diagrams.

If two knot diagrams are equivalent, the number of colourings using elements of
$S_m$ will be the same for all natural numbers $m$.
To show that two knot diagrams are not equivalent, it is sufficient to find a
value of $m$ such that the the invariant produces a different value for the two
knot diagrams.
Thus, the number of ways to colour both knot diagrams using elements of $S_m$
must be different.

It would be incorrect to simply count the number of ways that various values of
$S_m$ can be assigned to all $2n$ edges of the knot diagram.
To be an actual invariant, the assignments must satisfy two conditions.
The reason for that is that the values assigned to edges along an arc must be
the same.
Thus, if we have a crossing $X_{i,j,k,l}$, we label the values in $S_m$ that we
assign to each edge as $\sigma_i$, $\sigma_j$, $\sigma_k$, and $\sigma_l$
starting from the incoming lower edge and proceeding counterclockwise.
From this, the two relations that these four values have to satisfy are
constructed.

We know that the permutations must be equal along any arc.
Thus, the two edges of the top strand have the same permutations.

\papereq{Upper}{\sigma_j=\sigma_l}{}

The second criterion that our permutations must satisfy is as follows.
The signed product of the permutations assigned to the four edges around a
crossing must be the identity permutation.
To clarify, we start by picking an arbitrary sign convention, which in our case
is that inward pointing edges are given a positive sign.
Then, we move around a crossing in an arbitrary direction from an arbitrary
starting point, which in our case are counterclockwise from the edge labeled
$i$.
Inverting the permutations for the negative, and thus outward pointing, edges,
we construct an equation that positive crossings must satisfy.
Note that in a positive crossing, edges $j$ and $k$ point outwards.

\papereq{Positive}{\sigma_i\sigma_j^{\text-1}\sigma_k^{-1}\sigma_l=e}{}
\begin{paperwhere}
\papervar{e}{the identity permutation with $m$ elements}{}
\end{paperwhere}

For a negative crossing, we switch the signs of $j$ and $l$.

\papereq{Negative}{\sigma_i\sigma_j\sigma_k^{-1}\sigma_l^{\text-1}=e}{}

By putting together \eqUpper and \eqPositive, we get an equation for each
positive crossing.

\papereq{Right}{\sigma_k=\sigma_j\sigma_i\sigma_j^{\text-1}}{}

So we can find $\sigma_k$ by finding the permutation conjugation of $\sigma_i$
and $\sigma_j$.
We also get an equation for negative crossings.

\papereq{Left}{\sigma_k=\sigma_j^{\text-1}\sigma_i\sigma_j}{}

Thus, for negative crossings, $\sigma_k$ is the permutation conjugation of
$\sigma_i$ and $\sigma_j^{\text-1}$.

\pdfmsgex{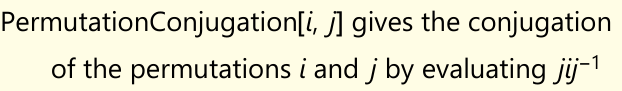}

We note that since we are taking the permutation conjugation of $\sigma_i$ with
$\sigma_j$, then $\sigma_i$ and $\sigma_k$ will have the same cycle lengths.
As any two edges across a crossing have the same cycle lengths, and since we can
follow the path of the knot by going through each crossing, one by one, never
changing cycle lengths, then all of the values for the edges must have the same
cycle lengths.
This gives us a lot more information.

Whereas before we would have had to map edges to $S_m$ and count the total
number of colourings, now they can be mapped to a subset of $S_m$.
If all the elements of this subset have the same cycle lengths, then the number
of colourings for each such subset can be counted independently.
Thus, instead of ending up with a single number as our invariant, we end up with
$P(m)$ different values, where $P$ is the partition function.
This is due to the fact that elements of $S_m$ have $P(m)$ different possible
cycle lengths.
Thus, to show that the knot diagrams are different, it suffices for any one of
these $P(m)$ values to differ.

As previously stated, we cannot simply find the number of ways to map edges to
elements of $S_m$.
For simplicity, we will assume that from now on we are dealing with a positive
crossing.
To do this for a negative crossing, it suffices to invert $\sigma_j$ before
computing its permutation conjugation with $\sigma_i$.

We need to find as many equations that relate the edges as possible using
\eqUpper, \eqPositive, and \eqNegative.
For knot diagrams with $n$ crossings we have to find equations relating the $2n$
edges.
Finding $n$ such equations is easy with \eqUpper as every crossing has an upper
arc whose two edges must have equal values in $S_m$.
Similarly, every crossing also gives us an equation from \eqRight or \eqLeft, so
we have a total of $2n$ equations.

We know that we cannot choose the permutations for all $2n$ edges independently
as there are equations relating them that will may not be satisfied if we were
to do so.
We also know choosing a permutation for just one edge will not be sufficient to
generate the remaining $2n-1$ edges.
Thus, there is some number of edges whose values can be used to generate the
rest of the $2n$ edges using \eqUpper, \eqRight, and \eqLeft.
We call such edges \textit{generators} and we need to find such a set of
edges.
Note that this set will not be unique, but we want to minimize its size while
still guaranteeing that we can generate all of the permutations.

For some perspective, when $n=10$, there are between 3 and 5 generators which
can determine the rest of the 20 permutations.
It is immediately clear that $n$ edges can be derived from the other $n$ by
using \eqUpper.
Thus, we are only interested in finding which of the other $n$ can generate all
of the permutations.
What we are effectively doing is writing out the indices from 1 to $2n$ and
striking out all those indices which can be determined from those that remain.
For each crossing, we examine the edges that make up the upper strand and remove
the one with the larger index from our list.
After we have done this, we have $n$ remaining indices.
We will call the set of these indices $S$.

We create a graph with vertices $V=\mathcal{P}(S)$, where $\mathcal{P}$ is the
power set function.
Note that $|V|=2^n$, which for $n=10$ is not too large.
For each crossing, \eqRight tells us that knowing $\sigma_j$ as well as either
$\sigma_i$ or $\sigma_k$, is sufficient to determine the other one, either
$\sigma_k$ or $\sigma_i$, respectively.
We represent this in the graph with a directed edge.
From every vertex with set $T\subset S$ where $T$ contains both $i$ and $j$, we
draw a directed edge to the vertex with set $T\cup\{k\}$.
Similarly, we draw an edge from each vertex whose set $T\subset S$ contains both
$j$ and $k$ to the vertex with set $T\cup\{i\}$.
An edge from vertex $v_1$ to vertex $v_2$ represents the fact that by knowing
the permutations for all of the edges with indices in the set $v_1$, we can
determine the permutations for all of the edges with indices in the set $v_2$ by
applying \eqRight or \eqLeft once.

We consider the vertices from which there is a path to $S$, the vertex
containing all of the indices.
If there is a path from a vertex $v$ to the vertex $S$, then we know that being
given the permutations for the edges whose indices are in $v$ is sufficient to
determine all of the permutations of the knot diagram.
Thus, to find our set of generators, we choose the vertex $v$ from which there
is a path from $v$ to $S$ such that $|v|$ is as small as possible.
To find the sequence of edges whose permutations can be determined, we examine
any path from $v$ to $S$ and choose that as our order.

\pdfmsgex{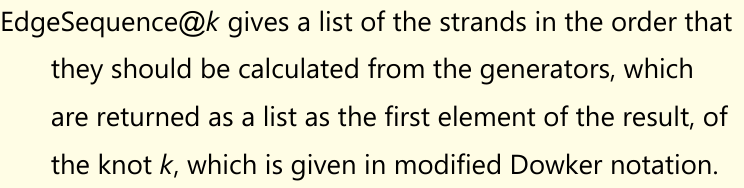}

Once we have found the generators, we need to start testing colourings.
We first break up $S_m$ into disjoint subsets by cycle lengths.
For each of these subsets $T$, we try assigning permutations in $T$ to each of
the generators in $v$ in every possible way.
There are $|T|^{|v|}$ ways of doing so.
This exponential is the reason for why we put so much effort into minimizing our
set of generators.

Once we set the generator permutations, we generate the permutations for the
rest of the edges, and then check that \eqRight or \eqLeft is satisfied for each
crossing.
If it is, then the colouring is valid, otherwise, it is not.

\pdfmsgex{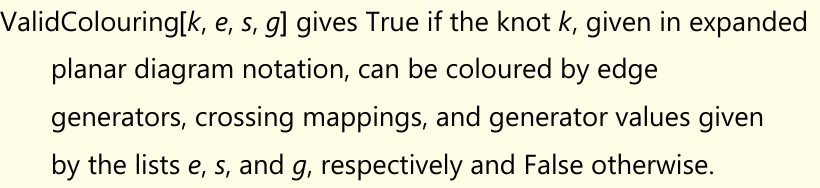}

We count the total number of valid colourings and our invariant becomes a list
of size $P(m)$ containing the number of valid colourings of the knot diagram
using elements of $S_m$.
Each element of the list corresponds to a different subset of $S_m$, where all
of the elements in each subset have the same cycle lengths.

After we had applied the Jones polynomial, we found two pairs of knot diagrams
that could be equivalent.
After setting $m$ to 5 and calculating the number of colourings for each of the
four knot diagrams in those two pairs, we found that the two knot diagrams in
our of the pairs produced different values, showing that they represent distinct
knots.
The two knots in the other pair produced the same list for the number of ways in
which they could be coloured using elements of $S_5$.
Thus, we have at least 250 distinct knot diagrams and potentially as many as
251.

\pdfmsgex{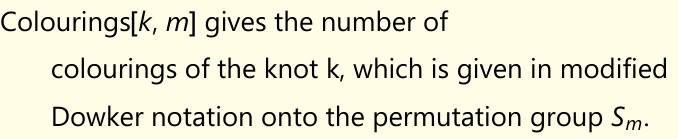}

We now try to see if the knot diagrams which have produced the same result for
each of our invariants are equivalent.
To show that the two knot diagrams in the pair are equivalent, we try adding a
positive kink into the knot diagrams in each of the 10 possible ways.
To add a kink, we insert a $k$ into the $k^\text{th}$ position in the MD code
and add 1 to all of the other values in the MD code that are greater than or
equal to $k$.

\pdfmsgex{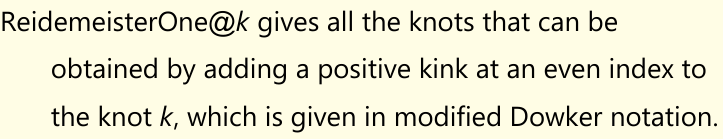}

Our two new 11-crossing knot diagrams are found to be equivalent under repeated
application of the third Reidemeister move (see \figMoves).

\pdfmsgex{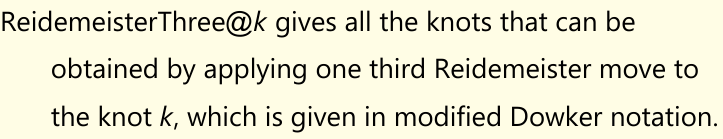}

Thus, the two knot diagrams that were producing the same values for each of our
invariants are equivalent.
This proves conclusively that we there are no fewer and no more than 250 knots
with up to 10 crossings.

\pdfmsgex{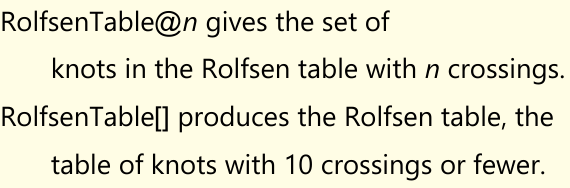}

\papersec{Knot Graphs}

In our calculation of the Rolfsen table, we mainly used three moves, the
2--pass, the third Reidemeister move, and the flype (see \figMoves).
We generated a graph of connections to find equivalent knot diagrams.
Now, we run our algorithms with a different set of knot diagrams.
We replace distinct prime alternating knot diagrams with all prime alternating
knot diagrams.
The difference between the sets is that the latter has diagrams that are
equivalent with respect to flyping.

From these knot diagrams, all of their non-alternating knot diagrams are
generated, each knot diagram is mapped to a vertex, these vertices are connected
with edges representing 2--passes, third Reidemeister moves, and flypes
(see \figMoves), and lastly all connected components that were found to be
reducible are removed.

\pdfmsgex{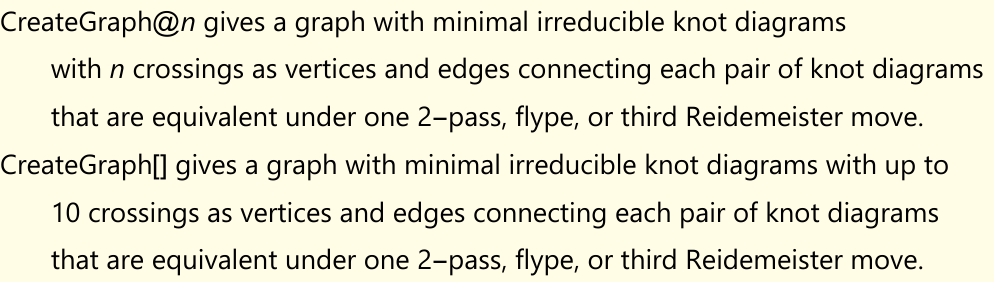}

The result is the full graph of irreducible knot diagrams and their connections.
This graph can be used for testing knot invariants.
Each invariant must produce the same result for each vertex in every connected
component.
Below is this graph, but with no labels.
The 2--passes, third Reidemeister moves, and flypes are represented by red,
green, and blue edges, respectively.

\vspace{-3em}
\begin{center}\scalebox{1.25}{\hspace{-0.1\columnwidth}\pdf{graph}}\end{center}
\vspace{-3em}

This graph is also available online at \url{http://tiny.cc/RolfsenTableGraph}.

\papersec{Utility Functions}

Here we include all the functions that are not mathematically interesting, but
merely serve as helper functions for those that are.
They are included here for completeness and in alphabetical order for ease of
access.

All of this code is also available online at
\url{http://tiny.cc/RolfsenTableCode}.

\pdfmsg{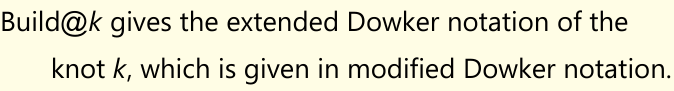}

\pdfmsg{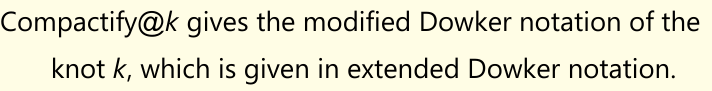}

\pdfmsg{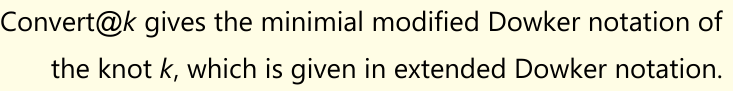}

\pdfmsg{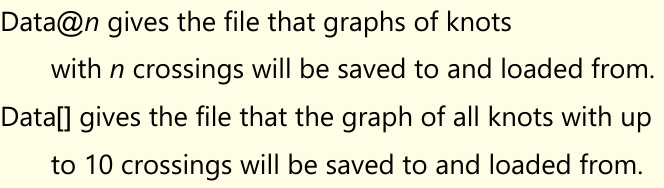}

\pdfmsg{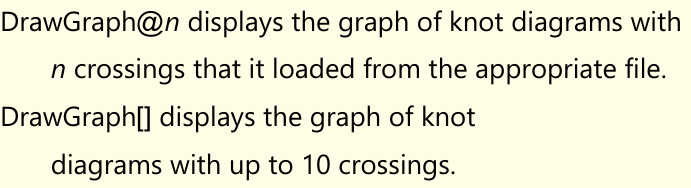}

\pdfmsg{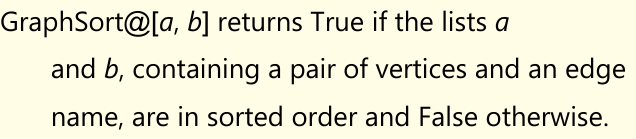}

\pdfmsg{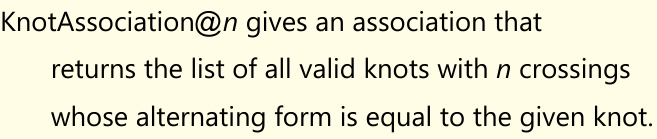}

\pdfmsg{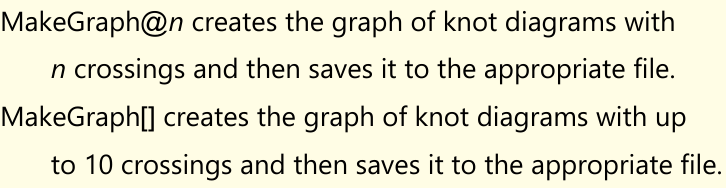}

\pdfmsg{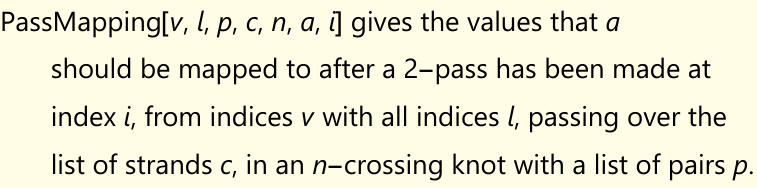}

\pdfmsg{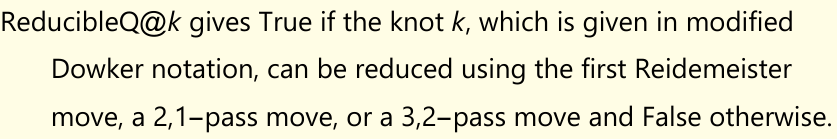}

\pdfmsg{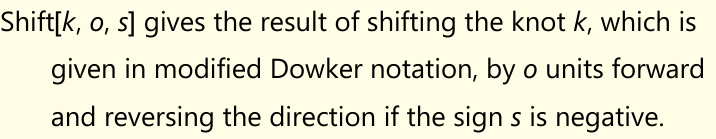}

\pdfmsg{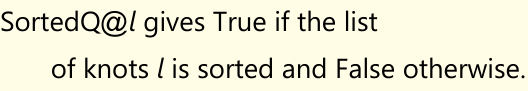}

\end{paper}

\fontsize{12pt}{12pt}\selectfont
\papersec{The Result}

\pdf{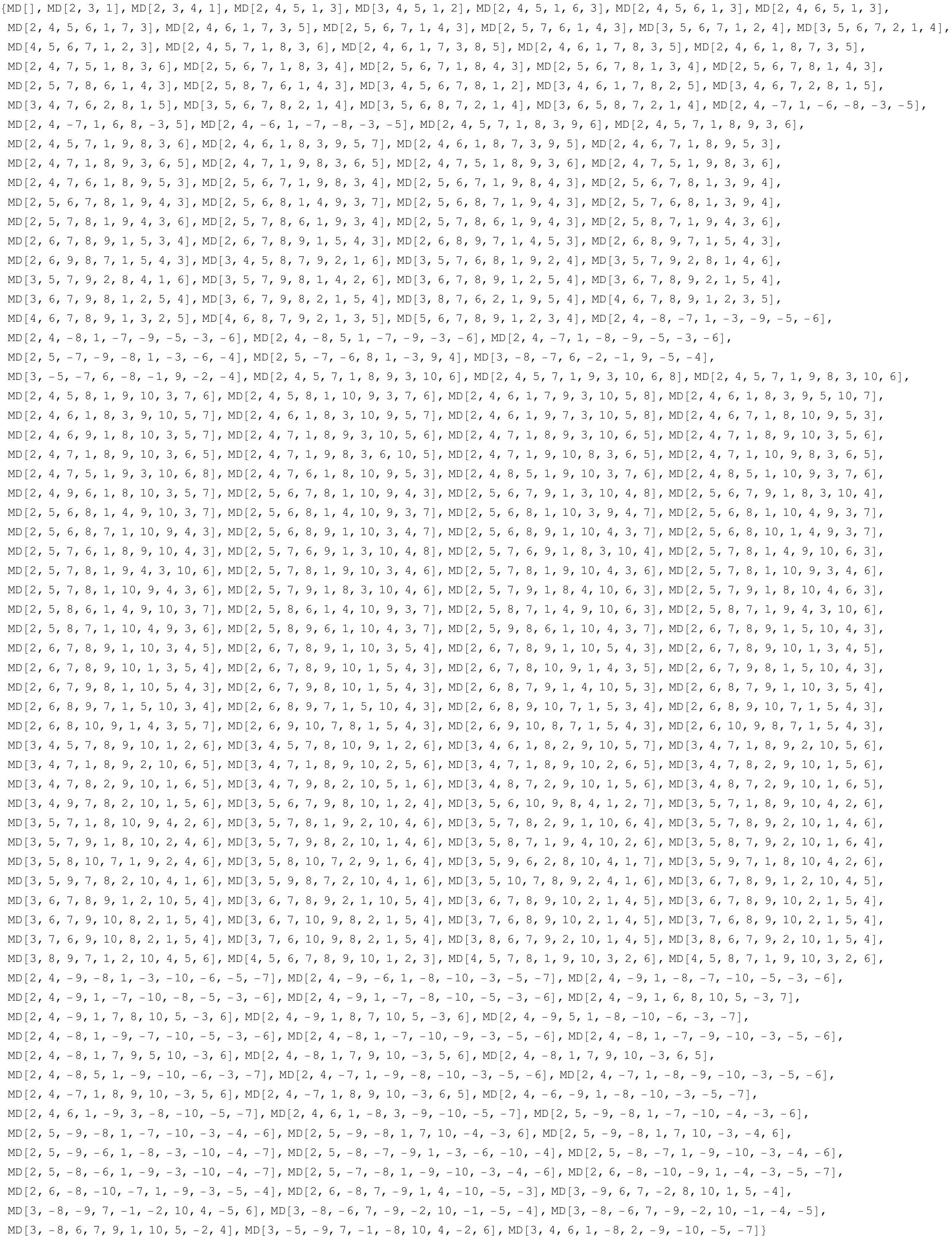}

This list can be found online at \url{http://tiny.cc/RolfsenTable}.

\begin{paper}

\papersec{Acknowledgements}

The author is grateful to Prof. Dror Bar-Natan for his neverending support and
assistance with this project.
Additionally, Prof. Bar-Natan originally suggested the topic of this project.

\papersec{References}

\end{paper}
\end{document}